\documentclass[english,smallcondensed,final]{svjour3}
\usepackage[T1]{fontenc}
\usepackage[latin9]{inputenc}
\usepackage{amssymb}
\usepackage{amsmath}
\usepackage{graphicx}
\usepackage{esint}
\usepackage{babel}
\def\d{{\mathbf{d}}}
\def\n{{\mathbf{n}}}
\def\u{{\mathbf{u}}}
\def\x{{\mathbf{x}}}
\def\y{{\mathbf{y}}}
\def\q{{\mathbf{q}}}

\smartqed
\providecommand{\tabularnewline}{\\}
\numberwithin{equation}{section}
\begin{document}
\title{A fast eikonal equation solver using the Schr\"odinger wave equation} 
\titlerunning{Linear embedding of the nonlinear eikonal equation}
\author {Karthik S. Gurumoorthy\thanks{This research work benefited from the support
of the AIRBUS Group Corporate Foundation Chair in Mathematics of Complex
Systems established in ICTS-TIFR} \and Adrian M. Peter \and Birmingham Hang Guan \and Anand Rangarajan}
\institute{Karthik S. Gurumoorthy \at
             International Center for Theoretical Sciences, Tata Institute of Fundamental Research, 
             TIFR Centre Building, Indian Institute of Science Campus, Bangalore, Karnataka, 560012, India\\
             \email{karthik.gurumoorthy@icts.res.in}
             \and
             Adrian M. Peter \at
             Department of Engineering Systems, Florida Institute of
             Technology, 150 W University Blvd., Melbourne, Florida, 32901,
             USA\\
             \email{apeter@fit.edu}
             \and
             Birmingham Hang Guan \at
             Department of Computer and Information Science and Engineering, University of Florida,
             E301 CSE Building, PO Box 116120, Gainesville, Florida, 32611, USA\\
             \email{hguan@cise.ufl.edu}
             \and
             Anand Rangarajan \at 
		  Department of Computer and Information Science and Engineering, University of Florida,
            E301 CSE Building, PO Box 116120, Gainesville, Florida, 32611, USA\\
            \email{anand@cise.ufl.edu}
}
\maketitle
\begin{abstract}
We use a Schr\"odinger wave equation formalism to solve the eikonal
equation. In our framework, a solution to the eikonal equation is obtained in
the limit as Planck's constant $\hbar$ (treated as a free parameter) tends to
zero of the solution to the corresponding linear Schr\"odinger equation. The
Schr\"odinger equation corresponding to the eikonal turns out to be a
\emph{generalized, screened Poisson equation}. Despite being linear, it does
not have a closed-form solution for arbitrary forcing functions. We present
two different techniques to solve the screened Poisson equation.  In the first
approach we use a standard perturbation analysis approach to derive a new
algorithm which is guaranteed to converge provided the forcing function is
bounded and positive. The perturbation technique requires a sequence of
discrete convolutions which can be performed in $O(N\log N)$ using the Fast
Fourier Transform (FFT) where $N$ is the number of grid points. In the second
method we discretize the linear Laplacian operator by the finite difference
method leading to a sparse linear system of equations which can be solved
using the plethora of sparse solvers.  The eikonal solution is recovered from
the exponent of the resultant scalar field. Our approach eliminates the need
to explicitly construct viscosity solutions as customary with direct solutions
to the eikonal. Since the linear equation is computed for a small but non-zero
$\hbar$, the obtained solution is an approximation. Though our solution
framework is applicable to the general class of eikonal problems, we detail
specifics for the popular vision applications of shape-from-shading, vessel
segmentation, and path planning.  

\end{abstract}

\section{Introduction}
\label{sec:Introduction}

The \emph{eikonal} (from the Greek word $\epsilon\iota\kappa o\nu$
or {}``image'') equation is traditionally encountered in the wave
and geometric optics literature where the principal concern is the
propagation of light rays in an inhomogeneous medium \cite{Chartier05}.
Its twin roots are in wave propagation theory and in geometric optics.
In wave propagation theory, it is obtained when the wave is approximated
using the Wentzel\textendash{}Kramers\textendash{}Brillouin (WKB)
approximation \cite{Paris69}. In geometric optics, it can be derived
from Huygens' principle \cite{Arnold89}. In the present day, the
eikonal equation has outgrown its humble optics origins and now finds
application in far flung areas such as electromagnetics \cite{Paris69},
robot motion path planning \cite{Canny99} and image analysis \cite{Osher02}.

The eikonal equation is a nonlinear, first order, partial
differential equation \cite{Whitham99} of the form 
\begin{equation}
\|\nabla S(\x)\|=f(\x),\, \x \in\Omega
\label{eikonalEq}
\end{equation}
subject to the boundary condition $S|_{\partial\Omega}=U(\x)$, where
$\Omega$ is an open subset of $\mathbb{R}^{D}$. The forcing function
$f(\x)$ is a positive valued function and $\nabla$ denotes
the gradient operator. In the special case where $f(\x)$ equals one
everywhere, the solution to the eikonal equation is the Euclidean
distance function \cite{Osher02}. Detailed discussions on the existence
and uniqueness of the solution can be found in \cite{Crandall92}.

While the eikonal equation is venerable and classical, it is only
in the last twenty years that we have seen the advent of numerical
methods aimed at solving this problem. To name a few are the pioneering
fast marching \cite{Osher88,Sethian96} and fast sweeping \cite{Zhao05} methods.
Algorithms based on discrete structures such as the well known Dijkstra
single source shortest path algorithm \cite{Cormen01} can also be
adapted to solve this problem. When we seek solutions on a discretized
spatial grid width $N$ points, the complexity of the fast marching
method is $O(N\log N)$ while that of the fast sweeping method for a single pass over the grid, is
$O(N)$ and therefore both of these efficient algorithms have seen
widespread use since their inception. The fast sweeping method is computationally nicer and easier to implement than the fast marching method, however the actual number of sweeps required for convergence depends on the problem at hand---experimentally it is observed that $2^D$ sweeps are required in $D$ dimensions. Recently, the ingenious work of Sapiro \emph{et al.} \cite{Yatziv06} provided an $O(N)$ implementation of the
fast marching method with a cleverly chosen untidy priority queue data structure. Typically, eikonal solvers
grow the solution from a set of $K$ seed points at which the solution
is known.

The eikonal equation can also be derived from a variational principle,
namely, Fermat's principle of least time which states that {}``Nature
always acts by the shortest paths'' \cite{Basdevant07}. From this
variational principle, the theoretical physics developmental sequence
proceeds as follows: The first order Hamilton's equations of motion
are derived using a Legendre transformation of the variational problem
wherein new momentum variables are introduced. Subsequently, a canonical
transformation converts the time varying momenta into constants of
the motion. The Hamilton-Jacobi equation emerges from the canonical
transformation \cite{Goldstein02}. In the Hamilton-Jacobi formalism
specialized to the eikonal problem, we seek a surface $S(X,t)$ such
that its increments are proportional to the speed of the light rays.
This is closely related to Huygens' principle and thus marks the \emph{rapprochement}
between geometric and wave optics \cite{Arnold89}. It is this nexus
that drives numerical analysis methods \cite{Osher88,Zhao05} (focused
on solving the eikonal equation) to base their solutions around the
Hamilton-Jacobi formalism.

So far, our development has followed that of classical physics. Since
the advent of quantum theory---specifically the Schr\"odinger wave equation---in
the 1920s, the close relationship between the Schr\"odinger and Hamilton-Jacobi
equations has been intensely studied \cite{Butterfield05}. Of particular
importance here is the quantum to classical transition as $\hbar\rightarrow0$
where the nonlinear Hamilton-Jacobi equation emerges from the phase
of the Schr\"odinger wave equation. This relationship has found very
few applications in the numerical analysis literature despite being
well known. In this paper, we leverage the important distinction between
the Schr\"odinger and Hamilton-Jacobi equations, namely, that the former
is \emph{linear} whereas the latter is not. We take advantage of the
linearity of the Schr\"odinger equation while exploiting its relationship
to Hamilton-Jacobi and derive computationally efficient solutions
to the eikonal equation.

A time-independent Schr\"odinger wave equation at the energy state \emph{E}
has the form $\hat{H}\phi(\x)=E \phi(\x)$
\cite{Griffiths04}, where $\phi(\x)$---the stationary state function---is the solution to the time-independent equation and $\hat{H}$ is the Hamiltonian operator. When the Hamilton-Jacobi scalar field $S^{\ast}$ is the exponent
of the stationary state function, specifically $\phi(\x)=\exp(\frac{-S^{\ast}(\x)}{\hbar})$,
and if $\phi(\x)$ satisfies the Schr\"odinger equation, we show that
as $\hbar\rightarrow0$, $S^{\ast}$ satisfies the Hamilton-Jacobi
equation. Note that in the above, a nonlinear Hamilton-Jacobi equation
is obtained in the limit as $\hbar\rightarrow0$ of a linear Schr\"odinger
equation which is novel from a numerical analysis perspective. Consequently,
instead of solving the Hamilton-Jacobi equation, one can solve its
Schr\"odinger counterpart (taking advantage of its linearity), and compute
an approximate $S^{\ast}$ for a suitably small value of $\hbar$.
This computational procedure is approximately equivalent to solving
the original Hamilton-Jacobi equation.

Since the efficient solution of a linear wave equation is the cornerstone
of our approach, we now briefly describe the actual computational
algorithm used. We derive the static Schr\"odinger equation for the
eikonal problem. The result is a generalized, \emph{screened Poisson}
equation \cite{Fetter03} whose solution is known at $K$ seed points.
This linear equation does not have a closed-form solution and therefore
we resort to a perturbation method \cite{Fernandez00} of solution---which
is related to the Born expansion \cite{Newton82}. The perturbation
method comprises a sequence of multiplications with a space-varying
forcing function followed by convolutions with a Green's
function (for the screened Poisson operator) which we solve using
an efficient $O(N\log N)$ fast Fourier transform (FFT)-based
technique \cite{Cooley65}. Perturbation analysis involves
a geometric series approximation for which we show convergence for
all bounded forcing functions independent of the value of $\hbar$. 

The intriguing characteristic of the perturbation approach is that it solves
the generalized screened Poisson without the explicit need to spatially
discretize the Laplacian operator. However, the downside of this method is
that it requires repeated convolution of the Green's function with the
solution from the previous iteration which can only be approximated by
discrete convolution. Hence the errors tend to accumulate with iteration. A
different route to solve the screened Poisson would be to approximate the
continuous Laplacian operator (say) by the method of finite differences and
use sparse, linear system solvers to obtain the solution. This approach leads
to many algorithm choices since there are myriad efficient sparse linear
solvers. We showcase the application of our linear discretized framework in
path planning, shape from shading and vessel segmentation.

The paper is organized as follows. In Section~\ref{sec:HJ formulation}, we
provide the Hamilton-Jacobi formulation for the eikonal equation as adopted by
the fast sweeping and fast marching methods. We restrict to the special case
of the eikonal equations involving constant forcing functions in
Section~\ref{sec:constforcingfunc} and derive its corresponding Schr\"odinger
wave equation. Section~\ref{sec:generalforcingfunc} considers the more general
version, where we derive and provide an efficient \emph{arbitrary precision}
FFT-based method for solving the Schr\"odinger equation using techniques form
perturbation theory. In Section~\ref{sec:discretizationApproach} we present a
second approach to solve the linear differential equation where by invoking a
finite difference approximation of the Laplacian operator we handle a sparse
linear system. We conclude in Section~\ref{Section:Discussion} by summarizing
our current work.

\section{Hamilton-Jacobi formulation for the eikonal equation}
\label{sec:HJ formulation}
\subsection{Fermat's principle of least time}
It is well known that the Hamilton-Jacobi equation formalism for the
eikonal equation can be obtained by considering a variational problem
based on Fermat's principle of least time \cite{Arnold89} which
in $2D$ is \begin{equation}
I[\q]=\int_{t_{o}}^{t_{1}}f(q_{1},q_{2},t)\sqrt{1+\dot{q}_{1}^{2}+\dot{q}_{2}^{2}}dt.\label{eq:I[q]}
\end{equation}

We take an idiosyncratic approach to the eikonal equation by considering
a different variational problem which is still very similar to Fermat's
least time principle. The advantage of this variational formulation
is that the corresponding Schr\"odinger wave equation can be easily
obtained.

Consider the following variational problem namely, 
\begin{equation}
I[\q]=\int_{t_{o}}^{t_{1}}\frac{1}{2}(\dot{q}_{1}^{2}+\dot{q}_{2}^{2})f^{2}(q_{1},q_{2})dt
\label{eq:I2[q]}
\end{equation}
 where the forcing term $f$ is assumed to be independent of time
and the Lagrangian $L$ is defined as \begin{equation}
L(q_{1},q_{2},\dot{q}_{1},\dot{q}_{2},t)\equiv\frac{1}{2}(\dot{q}_{1}^{2}+\dot{q}_{2}^{2})f^{2}(q_{1},q_{2}).\label{Lagrangian}\end{equation}
 Defining \begin{equation}
p_{i}\equiv\frac{\partial L}{\partial\dot{q}_{i}}=f^{2}(q_{1},q_{2})\dot{q}_{i}\label{eq:pi2}\end{equation}
 and applying the Legendre transformation \cite{Arnold89}, we can
obtain the Hamiltonian of the system in $2D$ as \begin{equation}
H(q_{1},q_{2},p_{1},p_{2},t)=\frac{1}{2}\frac{(p_{1}^{2}+p_{2}^{2})}{f^{2}(q_{1},q_{2})}.
\label{Hamiltonian}
\end{equation}

From a canonical transformation of the Hamiltonian \cite{Goldstein02},
we obtain the following Hamilton-Jacobi equation 
\begin{equation}
\frac{\partial S}{\partial t}+\frac{1}{2}\frac{\left(\frac{\partial S}{\partial q_{1}}\right)^{2}+\left(\frac{\partial S}{\partial q_{2}}\right)^{2}}{f^{2}(q_{1},q_{2})}=0\label{HJ}
\end{equation}
Since the Hamiltonian in (\ref{Hamiltonian}) is a constant \emph{independent}
of time, equation~(\ref{HJ}) can be simplified to the static Hamilton-Jacobi
equation. By separation of variables, we get 
\begin{equation}
S(q_{1},q_{2},t)=S^{\ast}(q_{1},q_{2})-Et
\label{eq:Sstar}
\end{equation}
where $E$ is the total energy of the system and $S^{\ast}(q_{1},q_{2})$
is called Hamilton's characteristic function \cite{Arnold89}. Observing
that $\frac{\partial S}{\partial q_{i}}=\frac{\partial S^{\ast}}{\partial q_{i}}$,
equation~(\ref{HJ}) can be rewritten as 
\begin{equation}
\frac{1}{2}\left[\left(\frac{\partial S^{\ast}}{\partial q_{1}}\right)^{2}+\left(\frac{\partial S^{\ast}}{\partial q_{2}}\right)^{2}\right]=Ef^{2}.\label{eq:HJstatic}
\end{equation}
 Choosing the energy $E$ to be $\frac{1}{2}$, we obtain 
\begin{equation}
\parallel\nabla S^{\ast}\parallel^{2}=f^{2}
\label{eq:gradS2}
\end{equation}
which is the original eikonal equation~(\ref{eikonalEq}). $S^{\ast}$
is the required Hamilton-Jacobi scalar field which is efficiently
obtained by the fast sweeping \cite{Zhao05} and fast marching methods
\cite{Osher88}.

\section{Eikonal equations with constant forcing functions}
\label{sec:constforcingfunc}
We begin the quantum formulation of the eikonal equation by considering its
special case where the forcing function is constant and equals $\tilde{f}$
everywhere.

\subsection{Deriving the Schr\"odinger wave equation}
 \label{sec:SchWaveEqconstf}

The time independent Schr\"odinger wave equation is given by \cite{Griffiths04} 
\begin{equation}
\hat{H}\phi(\x) = E \phi (\x)
\label{Schrodinger}
\end{equation}
where $\phi(\x)$ is the time-independent wave function and $\hat{H}$ is the Hamiltonian
operator obtained by first quantization where the momentum variables $p_{i}$ are replaced with the operator
$\frac{\hbar}{i}\frac{\partial}{\partial x_{i}}$. $E$ denotes the energy of the system. 

For this special case where the forcing functions is constant and equals
$\tilde{f}$ everywhere, the Hamiltonian of the system is given by (in 2D)
\begin{equation}
\label{eq:Hamiltonianconstf}
H(q_{1},q_{2},p_{1},p_{2},t)=\frac{1}{2}\frac{(p_{1}^{2}+p_{2}^{2})}{\tilde{f}^{2}(q_{1},q_{2})}.
\end{equation}
Its first quantization then yields the wave equation
\begin{equation}
\label{eq:TISchEqconstf}
-\frac{\hbar^{2}}{2\tilde{f}^2} \nabla^2 \phi = E \phi.
\end{equation}
When $E > 0$ we get oscillatory solution and when $E <0$ we get exponential
solutions in the sense of distributions.  In \cite{Gurumoorthy09,Sethi12} we
have shown that for the Euclidean distance function problem where $\tilde{f} =
1$, the exponential solution for $\phi$ obtained by setting $E = -\frac{1}{2}$
which is then used to recover $S^{\ast} $ using the relation
$\phi=\exp(\frac{-S^{\ast}}{\hbar})$, guarantees convergence of $S^{\ast}$ to
the true solution as $\hbar \rightarrow 0$. The work in \cite{Gurumoorthy09}
concerns with computing Euclidean distance functions only from point-sets, and
its extension of obtaining distance functions from curves is developed in
\cite{Sethi12}.  Following along similar lines, we propose to solve
(\ref{eq:TISchEqconstf}) at $E = -\frac{1}{2}$, for which $\phi$ satisfies the
differential equation
\begin{equation}
\label{eq:phiEquationconstf}
-\hbar^2 \nabla^2 \phi + \tilde{f}^2 \phi = 0.
\end{equation}

\subsection{Solving the Schr\"odinger wave equation}
\label{sec:computationconstf}
We now provide techniques for efficiently solving the Schr\"odinger equation
(\ref{eq:phiEquationconstf}). Note that we are interested in the computing the
solution only on the specified set of $N$ discrete grid locations.

Since the Laplacian operator $\nabla^{2}$ is negative definite, it follows
that the Hamiltonian operator $-\hbar^{2}\nabla^{2} + \tilde{f}^2$ is positive
definite for all values of $\hbar$ and hence the above equation does not have
a solution in the classical sense. Hence we look for a solution in the
distributional sense by considering the forced version of the equation, namely
\begin{equation}
-\hbar^{2}\nabla^{2}\phi+\tilde{f}^2\phi=\sum_{k=1}^{K}\delta(\x-\y_{k}).
\label{forcedPhiEquation}
\end{equation}
The points $\y_{k},k\in\{1,\ldots,K\}$ can be considered to be the
set of locations which encode initial knowledge about the scalar field
$S^{\ast}$, say for example $S^{\ast}(\y_k)=0,\forall \y_k,k\in\{1,\ldots,K\}$. 

For the forced equation (\ref{forcedPhiEquation}), closed-form solutions for
$\phi$ can be obtained in $1D$, $2D$ and $3D$ \cite{Gurumoorthy09} using the
Green's function approach \cite{Abramowitz64}. Since $S^{\ast}(\x)$ goes to
infinity for points at infinity, we can use Dirichlet boundary conditions
$\phi(\x)=0$ at the boundary of an unbounded domain. The form of the solution
for the Green's function $G$ is given by,

\noindent \textbf{1D:} In $1D$, the solution for $G$ \cite{Abramowitz64} is 
\begin{equation}
G(x)=\frac{1}{2\hbar\tilde{f}}\exp\left(\frac{-\tilde{f}|x|}{\hbar}\right).
\label{GreensFunction1D}
\end{equation}
 \textbf{2D:} In $2D$, the solution for $G$ \cite{Abramowitz64} is
 \begin{eqnarray}
G(\x) & = & \frac{1}{2\pi\hbar^{2}}K_{0}\left(\frac{\tilde{f}\|\x\|}{\hbar}\right)\label{GreensFunction2D}\\
 & \approx & \frac{\exp\left(\frac{-\tilde{f}\|\x\|}{\hbar}\right)}{2\hbar\sqrt{2\pi\hbar\tilde{f}\|\x\|}},\frac{\tilde{f}\|\x\|}{\hbar}\gg0.25\nonumber \end{eqnarray}
where $K_{0}$ is the modified Bessel function of the second kind.

\noindent \textbf{3D:} In $3D$, the solution for $G$ \cite{Abramowitz64} is 
\begin{equation}
G(\x)=\frac{1}{4\pi\hbar^{2}}\frac{\exp\left(\frac{-\tilde{f}\|\x\|}{\hbar}\right)}{\|\x\|}.
\label{GreensFunction3D}
\end{equation}

The solutions for $\phi$ can then be obtained by convolution 
\begin{equation}
\phi(\x)=\sum_{k=1}^{K}G(\x)\ast\delta(\x-\y_{k})=\sum_{k=1}^{K}G(\x-\y_{k})
\end{equation}
from which $S^{\ast}$ can be recovered using the relation (\ref{SfromPhi}). $S^{\ast}$ can explicitly be shown to converge to the the true solution $\tilde{f}r$, where $r=\min_{k}\|\x-\y_{k}\|$ as $\hbar\rightarrow0$ \cite{Gurumoorthy09}.

\subsubsection{Modified Green's function}
Based on the nature of the Green's function we would like to highlight on the following very important point. In the limiting case of $\hbar\rightarrow0$, 
\begin{equation}
\lim_{\hbar\rightarrow0}\frac{\exp\left\{\frac{-\tilde{f}\|\x\|}{\hbar}\right\} }{c\hbar^{d}\|\x\|^{p}}=0,\,\mathrm{for}\,\|\x\|\neq0\label{eq:gfbogacity}\end{equation}
for $c,d$ and $p$ being constants greater than zero and therefore we see that if we define
\begin{equation}
\tilde{G}(\x) = C \exp\left(\frac{-\tilde{f}\|\x\|}{\hbar}\right)
\end{equation} 
for some constant $C$, 
\begin{equation}
\lim_{\hbar\rightarrow0}|G(\x)-\tilde{G}(\x)|=0,\,\mathrm{for}\,\|\x\|\neq0
\label{eq:Gdiffbogacity}
\end{equation}
and furthermore, the convergence is \emph{uniform} for $\|\x\|$ away
from zero. Therefore, $\tilde{G}(\x)$ provides a very good approximation
for the actual Green's function as $\hbar \rightarrow 0$. For a fixed value of $\hbar$ and $\x$,
the difference between the Green's functions is $O\left(\frac{\exp\left(\frac{-\tilde{f}\|\x\|}{\hbar} \right )}{\hbar^2} \right)$
which is relatively insignificant for small values of $\hbar$ and for all $X \not= 0$. Moreover, using $\tilde{G}$ also avoids the singularity at the origin that $G$ has in the 2D and 3D case. The above observation motivates us to compute the solutions
for $\phi$ by convolving with $\tilde{G}$, namely
\begin{equation}
\phi(\x) = \sum_{k=1}^{K}\tilde{G}(\x) \ast \delta(\x-\y_{k})=\sum_{k=1}^{K}\tilde{G}(X-\y_k)
\end{equation}
 instead of the actual Green's function $G$ and recover $S^{\ast}$ using the relation (\ref{SfromPhi}), given by
\begin{equation}
S^{\ast}(\x)=-\hbar\log\left[\sum_{k=1}^{K}\exp\left(\frac{-\tilde{f}\|\x-\y_{k}\|}{\hbar}\right)\right] + \hbar \log(C).
\label{STilde}
\end{equation}
Since $\hbar \log(C)$ is a constant independent of $\x$ and converges to 0 as $\hbar \rightarrow 0$, it can ignored while computing $S^{\ast}$ at small values of $\hbar$--it is equivalent to setting $C$ to be 1. Hence the Schr\"odinger wave function for a constant force $\tilde{f}$ can be approximated by
\begin{equation} 
\label{eq:wavefuncconstf}
\phi(\x) = \sum_{k=1}^{K} \exp\left(\frac{-\tilde{f}\|X-\y_k\|}{\hbar}\right).
\end{equation}
It is worth emphasizing that the above defined wave function $\phi(\x)$ (\ref{eq:wavefuncconstf}), contains all the desirable properties that we need. Firstly, we notice that as $\hbar \rightarrow 0$, $\phi(\y_k) \rightarrow 1$ at the given point-set locations $\y_k$. Hence from (\ref{SfromPhi}) $S(\y_k) \rightarrow 0$ as $\hbar \rightarrow 0$ satisfying the initial conditions. Secondly as $\hbar\rightarrow0$, $\sum_{k=1}^{K}\exp\left(\frac{-\tilde{f}\|\x-\y_{k}\|}{\hbar}\right)$ can be approximated by $\exp\left(\frac{-\tilde{f}r}{\hbar}\right)$ where $r=\min_{k}\|\x-\y_{k}\|$. Hence $S^{\ast}(\x)$ $\approx-\hbar\log\exp\left(\frac{-\tilde{f}r}{\hbar}\right)=\tilde{f}r$, which is the true value. When $\tilde{f}=1$, we get the Euclidean distance function. Thirdly, $\phi$ can be easily computed using the fast Fourier transform as described under section (\ref{sec:computation}). Hence for computational purposes we consider the wave function defined in (\ref{eq:wavefuncconstf}).

\section{General eikonal equations}
\label{sec:generalforcingfunc}
Armed with the above set up, we can now solve the eikonal equation for arbitrary, positive-valued, bounded forcing functions $f$. We first show that even for general $f$, when $\phi$ satisfies the \emph{same} differential equation as in the case of constant forcing equation (replacing $\tilde{f}$ by $f$), namely
\begin{equation}
\label{eq:phiEquation}
-\hbar^2 \nabla^2 \phi + f^2 \phi = 0,
\end{equation}
and is related to $S^{\ast}$ by $\phi=\exp(\frac{-S^{\ast}}{\hbar})$, $S^{\ast}$ asymptotically satisfies the eikonal equation (\ref{eikonalEq}) as $\hbar\rightarrow0$. We show this for the $2D$ case but the generalization
to higher dimensions is straightforward.

When $\phi(x_{1},x_{2})=\exp(\frac{-S^{\ast}(x_{1},x_{2})}{\hbar})$,
the first partials of $\phi$ are 
\begin{equation}
\frac{\partial\phi}{\partial x_{1}}=-\frac{1}{\hbar}\exp\left(\frac{-S^{\ast}}{\hbar}\right)\frac{\partial S^{\ast}}{\partial x_{1}},\frac{\partial\phi}{\partial x_{2}}=-\frac{1}{\hbar}\exp\left(\frac{-S^{\ast}}{\hbar}\right)\frac{\partial S^{\ast}}{\partial x_{2}}.
\end{equation}
 The second partials required for the Laplacian are 
\begin{eqnarray}
\frac{\partial^{2}\phi}{\partial x_{1}^{2}}=\frac{1}{\hbar^{2}}\exp\left(\frac{-S^{\ast}}{\hbar}\right)\left(\frac{\partial S^{\ast}}{\partial x_{1}}\right)^{2}-\frac{1}{\hbar}\exp\left(\frac{-S^{\ast}}{\hbar}\right)\frac{\partial^{2}S^{\ast}}{\partial x_{1}^{2}},\nonumber \\
\frac{\partial^{2}\phi}{\partial x_{2}^{2}}=\frac{1}{\hbar^{2}}\exp\left(\frac{-S^{\ast}}{\hbar}\right)\left(\frac{\partial S^{\ast}}{\partial x_{2}}\right)^{2}-\frac{1}{\hbar}\exp\left(\frac{-S^{\ast}}{\hbar}\right)\frac{\partial^{2}S^{\ast}}{\partial x_{2}^{2}}.
\end{eqnarray}
 From this, equation~(\ref{eq:phiEquation}) can be rewritten as 
\begin{equation}
\left(\frac{\partial S^{\ast}}{\partial x_{1}}\right)^{2}+\left(\frac{\partial S^{\ast}}{\partial x_{2}}\right)^{2}-\hbar\left(\frac{\partial^{2}S^{\ast}}{\partial x_{1}^{2}}+\frac{\partial^{2}S^{\ast}}{\partial x_{2}^{2}}\right)=f^{2}\end{equation}
 which in simplified form is 
\begin{equation}
\|\nabla S^{\ast}\|^{2}-\hbar\nabla^{2}S^{\ast}=f^{2}.
\label{Sequation}
\end{equation}
The additional $\hbar\nabla^{2}S^{\ast}$ term {[}relative
to (\ref{eikonalEq}){]} is referred to as the \emph{viscosity term}
\cite{Crandall92,Osher88} which emerges naturally from the Schr\"odinger
equation derivation---an intriguing result which differs from direct solutions of the non-linear eikonal that \emph{artificially} incorporate viscosity terms. Since $|\nabla^{2}S^{\ast}|$
is bounded, as $\hbar\rightarrow0$, (\ref{Sequation}) tends to $\|\nabla S^{\ast}\|^{2}=f^{2}$
which is the original eikonal equation~(\ref{eikonalEq}). This relationship
motivates us to solve the linear Schr\"odinger equation~(\ref{eq:phiEquation})
instead of the non-linear eikonal equation and then compute the scalar
field $S^{\ast}$ via 
\begin{equation}
S^{\ast}(\x)=-\hbar\log\phi(\x).
\label{SfromPhi}
\end{equation}

\section{Perturbation theory approach to solver the linear system}
 Since the linear system (\ref{eq:phiEquation}) (and its forced version) doesn't have a closed-form solution for non-constant forcing functions, one approach to solve it is using \emph{perturbation} theory \cite{Fernandez00}. Assuming that $f$ is close to a constant non-zero forcing function $\tilde{f}$, equation~(\ref{eq:phiEquation}) can be rewritten as 
\begin{equation}
(-\hbar^{2}\nabla^{2}+\tilde{f}^{2})\left[1+(-\hbar^{2}\nabla^{2}+\tilde{f}^{2})^{-1}\circ(f^{2}-\tilde{f}^{2})\right]\phi=0.
\end{equation}
 Now, defining the operator 
 \begin{equation}
L\equiv(-\hbar^{2}\nabla^{2}+\tilde{f}^{2})^{-1}\circ(f^{2}-\tilde{f}^{2})
\label{Ldef}
\end{equation}
 and the function $\phi_{0}\equiv(1+L)$, we see that $\phi_{0}$ satisfies 
\begin{equation}
(-\hbar^{2}\nabla^{2}+\tilde{f}^{2})\phi_{0}=0
\label{phi_0Equation}
\end{equation}
 and 
\begin{equation}
\phi=(1+L)^{-1}\phi_{0}.
\label{eq:phiandphi_0Equation}
\end{equation}
Notice that in the differential equation for $\phi_0$ (\ref{phi_0Equation}), the forcing function is constant and equals $\tilde{f}$ everywhere. Hence $\phi_0$ behaves like the wave function corresponding to the constant forcing function $\tilde{f}$---described under section~(\ref{sec:computationconstf}) and can be approximated by 
\begin{equation}
\label{eq:phi_0Solution}
\phi_0(\x) = \sum_{k=1}^{K}\exp\left(\frac{-\tilde{f}\|\x-\y_{k}\|}{\hbar}\right).
\end{equation}
We solve for $\phi$ in (\ref{eq:phiandphi_0Equation}) using a geometric series approximation for $(1+L)^{-1}$. Firstly, observe that the approximate solution for $\phi_0$ in (\ref{eq:phi_0Solution}) is a square-integrable function which is necessary for the subsequent steps.

Let $\mathcal{H}$ denote the space of square integrable functions on $\mathbb{R}^{D}$, i.e, $g\in\mathcal{H}$ if and only if $\int g^{2}d\mu<\infty$.
The function norm $\|g\|$ for a function $g\in\mathcal{H}$ is given by $\|g\|^{2}=\int g^{2}d\mu$ where $\mu$ is the Lebesgue measure on $\mathbb{R}^{D}$. Let $\mathcal{B} = \{g\in\mathcal{H}:\|g\|\leq1\}$ denote a closed unit ball in the Hilbert space $\mathcal{H}$ and $c_{0} \equiv \|L\|_{op}$ be the operator norm defined as 
\begin{equation}
c_{0}=\sup\{\|Lg\|,\, \forall g \in \mathcal{B}\}.
\label{NormDefinition}
\end{equation}
If $c_{0}<1$, we can approximate $(1+L)^{-1}$ using the first few $T+1$ terms of the geometric series to get 
\begin{equation}
(1+L)^{-1}\approx1-L+L^{2}-L^{3}+\ldots+(-1)^{T}L^{T}
\label{invApprox}
\end{equation}
 where the operator norm of the difference can be bounded by 
\begin{equation}
\|(1+L)^{-1}-\sum_{i=0}^{T}(-1)^{i}L^{i}\|_{op}\leq\sum_{i=T+1}^{\infty}\|L^{i}\|_{op}\leq\sum_{i=T+1}^{\infty}c_{0}^{i}=\frac{c_{0}^{T+1}}{1-c_{0}}
\label{OperatorNormDiff}
\end{equation}
which converges to $0$ \emph{exponentially} in $T$. We would like
to point out that the above geometric series approximation is similar
to a Born expansion used in scattering theory \cite{Newton82}. We
now derive an upper bound for $c_{0}$.

Let $L=A_{1}\circ A_{2}$ where $A_{1}\equiv(-\hbar^{2}\nabla^{2}+\tilde{f}^{2})^{-1}$
and $A_{2}\equiv f^{2}-\tilde{f}^{2}$. We now provide an upper bound
for $\|A_{1}\|_{op}$. For a given $g\in\mathcal{B}$, let $z = A_1(g)$ where $z$ satisfies the relation $(-\hbar^{2}\nabla^{2}+\tilde{f}^{2})z=g$ with vanishing Dirichlet boundary conditions at $\infty$. Then
\begin{eqnarray}
\|(-\hbar^{2}\nabla^{2}+\tilde{f}^{2})z\|^{2} & = & \|-\hbar^{2}\nabla^{2}z\|^{2}+\|\tilde{f}^{2}z\|^{2}+2\hbar^{2}\tilde{f}^{2}\langle-\nabla^{2}z,z\rangle\nonumber \\
 & = & \|g\|^{2}\leq1.\label{eq:functionNorm}
\end{eqnarray}
Using the identity $\nabla.(z\nabla z)=z\nabla^{2}z+|\nabla z|^{2}$ we write
\begin{equation}
\langle-\nabla^{2}z,z\rangle=-\int z\nabla^{2}zd\mu=-\int\nabla.(z\nabla z)d\mu+\int|\nabla z|^{2}d\mu.
\end{equation}
Divergence theorem states that $-\int\nabla.(z\nabla z)d\mu=0$ and hence 
\begin{equation}
\langle-\nabla^{2}z,z\rangle=\int|\nabla z|^{2}d\mu\geq0.
\end{equation}
Using the above relation in (\ref{eq:functionNorm}) we find $ \|z\|=\|A_{1}(g)\| \leq\frac{1}{\tilde{f}^{2}},\forall g \in \mathcal{B}$ implying that
\begin{equation}
\|A_{1}\|_{op}\leq\frac{1}{\tilde{f}^{2}}.
\label{boundOperator1}
\end{equation}
Furthermore, as for any $g\in\mathcal{B}$ 
\begin{equation}
\|(f^{2}-\tilde{f}^{2})g\|^{2}=\int(f^{2}-\tilde{f}^{2})^{2}g^{2}d\mu\leq \left(\sup\{|f^{2}-\tilde{f}^{2}|\}\right)^2
\end{equation}
 we get
\begin{equation}
\|A_{2}\|_{op}\leq \sup\{|f^{2}-\tilde{f}^{2}|\}.
\label{boundOperator2}
\end{equation}
 Since $\|L\|_{op}\leq\|A_{1}\|_{op}\|A_{2}\|_{op}$, from equations~(\ref{boundOperator1})
and (\ref{boundOperator2}) we can deduce that
\begin{equation}
c_{0}=\|L\|_{op}\leq\frac{\sup\{|f^{2}-\tilde{f}^{2}|\}}{\tilde{f}^{2}}.
\label{OperatorNormUpperBound}
\end{equation}
It is worth commenting that the bound for $c_{0}$ is \emph{independent}
of $\hbar$. So, if we guarantee that $\frac{\sup\{|f^{2}-\tilde{f}^{2}|\}}{\tilde{f}^{2}}<1$,
the geometric series approximation for $(1+L)^{-1}$ (\ref{invApprox})
converges for \emph{all} values of $\hbar$.

\subsection{Deriving a bound for convergence of the perturbation series}
Interestingly, for any positive, upper bounded forcing function $f$
bounded away from zero, i.e $f(\x)>\epsilon$ for some $\epsilon>0$%
\footnote{If $f(\x)=0$, then the velocity $v(\x)=\frac{1}{f(\x)}$ becomes $\infty$
at $\x$. Hence it is reasonable to assume $f(\x)>0$.%
}, by defining $\tilde{f}=\sup\{f(\x)\}$, we observe that $|f^{2}-\tilde{f}^{2}|<\tilde{f}^{2}$.
From equation~(\ref{OperatorNormUpperBound}), we immediately see
that $c_{0}<1$. This proves the existence of $\tilde{f}$ for which
the geometric series approximation (\ref{invApprox}) is \emph{always}
guaranteed to converge for any positive bounded forcing function $f$
bounded away from zero. The choice of $\tilde{f}$ can then be made
prudently by defining it to be the value that \emph{minimizes} 
\begin{equation}
F(\tilde{f})=\frac{\sup\{|f^{2}-\tilde{f}^{2}|\}}{\tilde{f}^{2}}.
\end{equation}
 This in turn minimizes the operator norm $c_{0}$, thereby providing
a better geometric series approximation for the inverse (\ref{invApprox}).

Let $f_{min}=\inf\{f(\x)\}$ and let $f_{max}=\sup\{f(\x)\}$. We now
show that $F(\tilde{f})$ attains its minimum at \begin{equation}
\tilde{f}=\nu=\sqrt{\frac{f_{min}^{2}+f_{max}^{2}}{2}}.\label{fTildeValue}\end{equation}
case (i): If $\tilde{f}<\nu$, then $\sup\{|f^{2}-\tilde{f}^{2}|\}=f_{max}^{2}-\tilde{f}^{2}$.
Clearly, \begin{equation}
\frac{f_{max}^{2}-\tilde{f}^{2}}{\tilde{f}^{2}}>\frac{f_{max}^{2}-\nu^{2}}{\nu^{2}}.\end{equation}
case (ii): If $\tilde{f}>\nu$, then $\sup\{|f^{2}-\tilde{f}^{2}|\}=\tilde{f}^{2}-f_{min}^{2}$.
It follows that \begin{equation}
\frac{\tilde{f}^{2}-f_{min}^{2}}{\tilde{f}^{2}}=1-\frac{f_{min}^{2}}{\tilde{f}^{2}}>1-\frac{f_{min}^{2}}{\nu^{2}}.\end{equation}
 We therefore see that $\tilde{f}=\nu=\sqrt{\frac{f_{min}^{2}+f_{max}^{2}}{2}}$
is the optimal value.

Using the above approximation for $(1+L)^{-1}$ (\ref{invApprox})
and the definition of $L$ from (\ref{Ldef}) we obtain the solution
for $\phi$ as 
\begin{equation}
\phi=\phi_{0}-\phi_{1}+\phi_{2}-\phi_{3}+\ldots+(-1)^{T}\phi_{T}
\label{eq:phiSolution}
\end{equation}
 where $\phi_{i}$ satisfies the recurrence relation 
\begin{equation}
(-\hbar^{2}\nabla^{2}+\tilde{f}^{2})\phi_{i}=(f^{2}-\tilde{f}^{2})\phi_{i-1},\hspace{2pt}\forall i\in\{1,2,\ldots,T\}.
\label{phiRecRelation}
\end{equation}
Observe that (\ref{phiRecRelation}) is an inhomogeneous, screened Poisson equation with a constant forcing function $\tilde{f}$. Following a Green's function approach \cite{Abramowitz64}, each $\phi_i$ can be obtained by convolution
\begin{equation}
\label{eq:phi_i_Solution}
\phi_i = G\ast\left[(f^{2}-\tilde{f}^{2})\phi_{i-1}\right]
\end{equation}
where $G$ is given by equations (\ref{GreensFunction1D}), (\ref{GreensFunction2D}) or (\ref{GreensFunction3D}) depending upon the spatial dimension.

Once the $\phi_{i}$'s are computed, the wave function $\phi$ can then
be determined using the approximation (\ref{eq:phiSolution}). The solution
for the eikonal equation can be recovered using the relation (\ref{SfromPhi}). Notice that if $f = \tilde{f}$ everywhere, then all $\phi_i$'s except $\phi_0$ is identically equal to zero and we get $\phi = \phi_0$ as described under section (\ref{sec:constforcingfunc}).
\subsection{Efficient computation of the wave function}
\label{sec:computation}
In this section, we provide numerical techniques for efficiently computing the wave function $\phi$. Recall that we are interested in solving the eikonal equation only at the given $N$ discrete grid locations. 
Consider the solution for $\phi_0$ given in (\ref{eq:phi_0Solution}). In order to obtain the desired solution for $\phi_0$ computationally, we must replace the $\delta$ function by the Kronecker delta function
\begin{equation}
\label{eq:deltakron}
\delta_{kron}(\x) = \left\{ 
	\begin{array}{ll}
         1 & \mbox{if $\x = \y_k$};\\
        0 & \mbox{otherwise}
        \end{array} 
        \right. 
\end{equation}
that takes $1$ at the point-set locations ($\{\y_k\}$) and $0$ at other grid locations. Then $\phi_0$ can be \emph{exactly} computed at the grid locations by the discrete convolution of $\tilde{G}$ (setting $C=1$) with the Kronecker-delta function.

To compute $\phi_i$, we replace each of the convolutions in (\ref{eq:phi_i_Solution}) with the discrete convolution between the functions computed at the $N$ grid locations. By the convolution theorem \cite{Bracewell99}, a
discrete convolution can be obtained as the inverse Fourier transform
of the product of two individual transforms which for two $O(N)$
sequences can be performed in $O(N\log N)$ time \cite{Cooley65}.
Thus, the values of each $\phi_{i}$ at the $N$ grid locations can
be efficiently computed in $O(N\log N)$ making use of the values
of $\phi_{i-1}$ determined at the earlier step. Thus, the overall
time complexity to compute the approximate $\phi$ using the first
few $T+1$ terms is then $O(TN\log N)$. Taking the logarithm of $\phi$
then provides an approximate solution to the eikonal equation. The
algorithm is summarized in Table~\ref{table:algorithm}.

\begin{center}
\begin{table}
\caption{Algorithm for the approximate solution of the eikonal equation}
\label{table:algorithm} 
\centering{}
\begin{tabular}{cl}
\hline 
1.  & Compute the function $\tilde{G}(\x)=\exp\left(\frac{-\tilde{f}\|\x\|}{\hbar}\right)$
at the grid locations.\tabularnewline
2.  & Define the function $\delta_{kron}(\x)$ which takes the value $1$ at the point-set
locations \tabularnewline
 & and $0$ at other grid locations.\tabularnewline
3.  & Compute the FFT of $\tilde{G}$ and $\delta_{kron}$, namely $\tilde{G}_{FFT}(\u)$
and $\delta_{FFT}(\u)$ respectively.\tabularnewline
4.  & Compute the function $H(\u)=\tilde{G}_{FFT}(\u)\delta_{FFT}(\u)$.\tabularnewline
5.  & Compute the inverse FFT of $H$ to obtain $\phi_0(\x)$
at the grid locations.\tabularnewline
6.  & Initialize $\phi(\x)$ to $\phi_{0}(\x)$.\tabularnewline
7. &Consider the Green's function $G$ corresponding to the spatial dimension \tabularnewline 
   & and compute its FFT, namely $G_{FFT}(\u)$.\tabularnewline 
8.  & For $i\hspace{2pt}=\hspace{2pt}1$ to $T$ do \tabularnewline
9.  & \hspace{0.5in}Define $\psi_i(\x)=\left[f^{2}(\x)-\tilde{f}^{2}\right]\phi_{i-1}(\x)$.\tabularnewline
10.  & \hspace{0.5in}Compute the FFT of $\psi_i$ namely $\Psi(\u)$. \tabularnewline
11.  & \hspace{0.5in}Compute the function $H(\u)=G_{FFT}(\u) \Psi_i(\u)$.\tabularnewline
12.  & \hspace{0.5in}Compute the inverse FFT of $H$ and multiply it with
the grid width\tabularnewline
 & \hspace{0.5in}area/volume to compute $\phi_{i}(\x)$ at the
grid locations.\tabularnewline
13.  & \hspace{0.5in}Update $\phi(\x)=\phi(\x)+(-1)^{i}\phi_{i}(\x)$.\tabularnewline
14.  & End \tabularnewline
15.  & Take the logarithm of $\phi(\x)$ and multiply it by $(-\hbar)$
to get \tabularnewline
 & the approximate solution for the eikonal equation at the grid locations.\tabularnewline
\hline
\end{tabular}
\end{table}
\par\end{center}
We would like to emphasize that the number of terms ($T$) used in
the geometric series approximation of $(1+L)^{-1}$ (\ref{invApprox})
is \emph{independent} of $N$. Using more terms only improves the
approximation of this truncated geometric series as shown in the experimental section. 
From equation~(\ref{OperatorNormDiff}),
it is evident that the error incurred due to this approximation converges
to zero \emph{exponentially} in $T$ and hence even with a small value
of $T$, we should be able to achieve good accuracy.

\subsection{Numerical issues}
In principle, we should be able to apply our technique at very small
values of $\hbar$ and obtain highly accurate results. But we noticed
that a na\"{i}ve double precision-based implementation tends to deteriorate
for $\hbar$ values very close to zero. This is due to the fact that
at small values of $\hbar$ (and also at large values of $\tilde{f}$),
$\exp\left(\frac{-\tilde{f}\|\x\|}{\hbar}\right)$ drops off very quickly
and hence for grid locations which are far away from the point-set,
the convolution done using FFT may not be accurate. To this end, we
turned to the GNU MPFR multiple-precision arithmetic library which
provides arbitrary precision arithmetic with correct rounding \cite{Fousse07}. 
MPFR is based on the GNU multiple-precision library (GMP) \cite{GMP}. It enabled us
to run our technique at very small values of $\hbar$ giving highly
accurate results. We corroborate our claim and demonstrate the usefulness
of our method with the set of experiments described in the subsequent
section.

\subsection{Experimental verification of the perturbation approach}
\label{sec:Experiments}
In this section, we demonstrate the usefulness of our perturbation approach by
computing the approximate solution to the general eikonal equation~(\ref{eikonalEq})
over a $2D$ grid.
\subsubsection{Comparison with the true solution}
\noindent \textbf{Example 1:} In this example, we solve the eikonal
equation for the scenario where the exact solution is known \emph{a
priori} at the grid locations. The exact solution is $R(x,y)=|e^{\sqrt{x^{2}+y^{2}}}-1|$. The boundary condition is $R(x,y)=0$ at the point source located at $(x_{0},y_{0})=(0,0)$. The forcing function---the absolute gradient
$|\nabla R|$---is $f(x,y)=|\nabla R|=e^{\sqrt{x^{2}+y^{2}}}$ specified on a $2D$ grid consisting of points between $(-0.125,-0.125)$
and $(0.125,0.125)$ with a grid width of $\frac{1}{2^{10}}$. We
ran the Schr\"odinger for $6$ iterations at $\hbar=0.006$ and the
fast sweeping for $15$ iterations sufficient enough for both the
methods to converge. The percentage error is calculated according to
\begin{equation}
\Lambda=\frac{100}{N}\sum_{i=1}^{N}\frac{\Delta_{i}}{R_{i}},\label{eqn:percentError}
\end{equation}
 where $R_{i}$ and $\Delta_{i}$ are respectively the actual value and the absolute difference of the computed and actual
value at the $i^{th}$ grid point. The maximum difference between the true and approximate solution for different iterations is summarized in the table~(\ref{table:eikonalexp1}).

\noindent \begin{center}
\begin{table}[ht!] 
\centering{}
\caption{Percentage error for the Schr\"odinger method for different iterations}
\label{table:eikonalexp1}
\begin{tabular}{|c|c|c|}
\hline 
Iter  & \% error & max diff \tabularnewline
\hline 
1 & 2.081042 & 0.002651\tabularnewline
\hline 
2 & 1.514745 & 0.002140 \tabularnewline
\hline 
3 & 1.390552 & 0.002142\tabularnewline
\hline 
4 & 1.363256 & 0.002128 \tabularnewline
\hline 
5 & 1.357894 & 0.002128 \tabularnewline
\hline 
6 & 1.356898 & 0.002128 \tabularnewline
\hline 
\end{tabular}
\end{table}
\end{center}
The fast sweeping gave a percentage error of \emph{$1.135\%$}. We believe that the error incurred in our Schr\"odinger approach can be further reduced by decreasing $\hbar$ but at the expense of more computational power requiring higher precision floating point arithmetic.

\begin{figure}[ht!]
\begin{centering}
\begin{tabular}{ccc}
\includegraphics[width=0.3\textwidth]{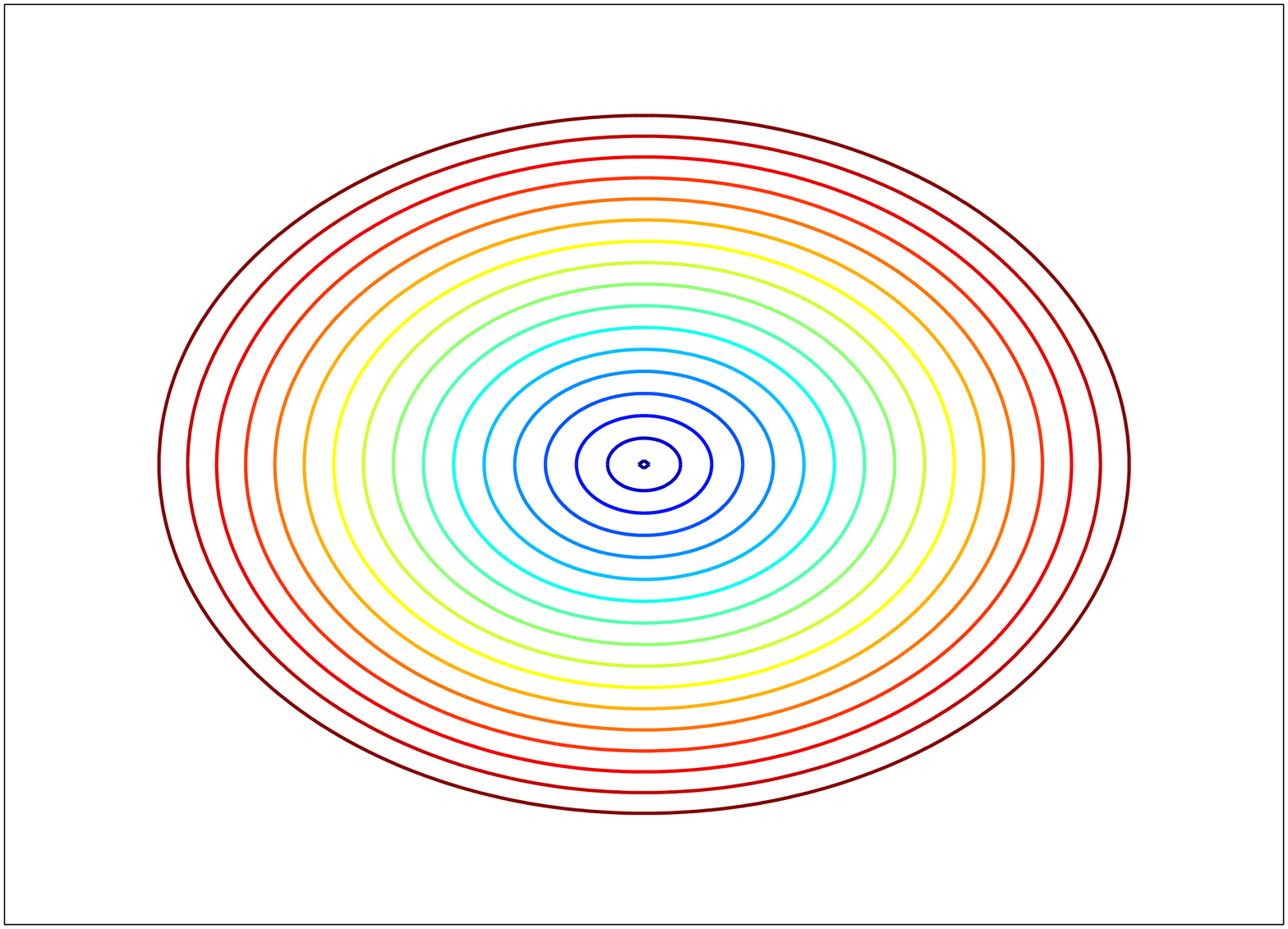}&
\includegraphics[width=0.3\textwidth]{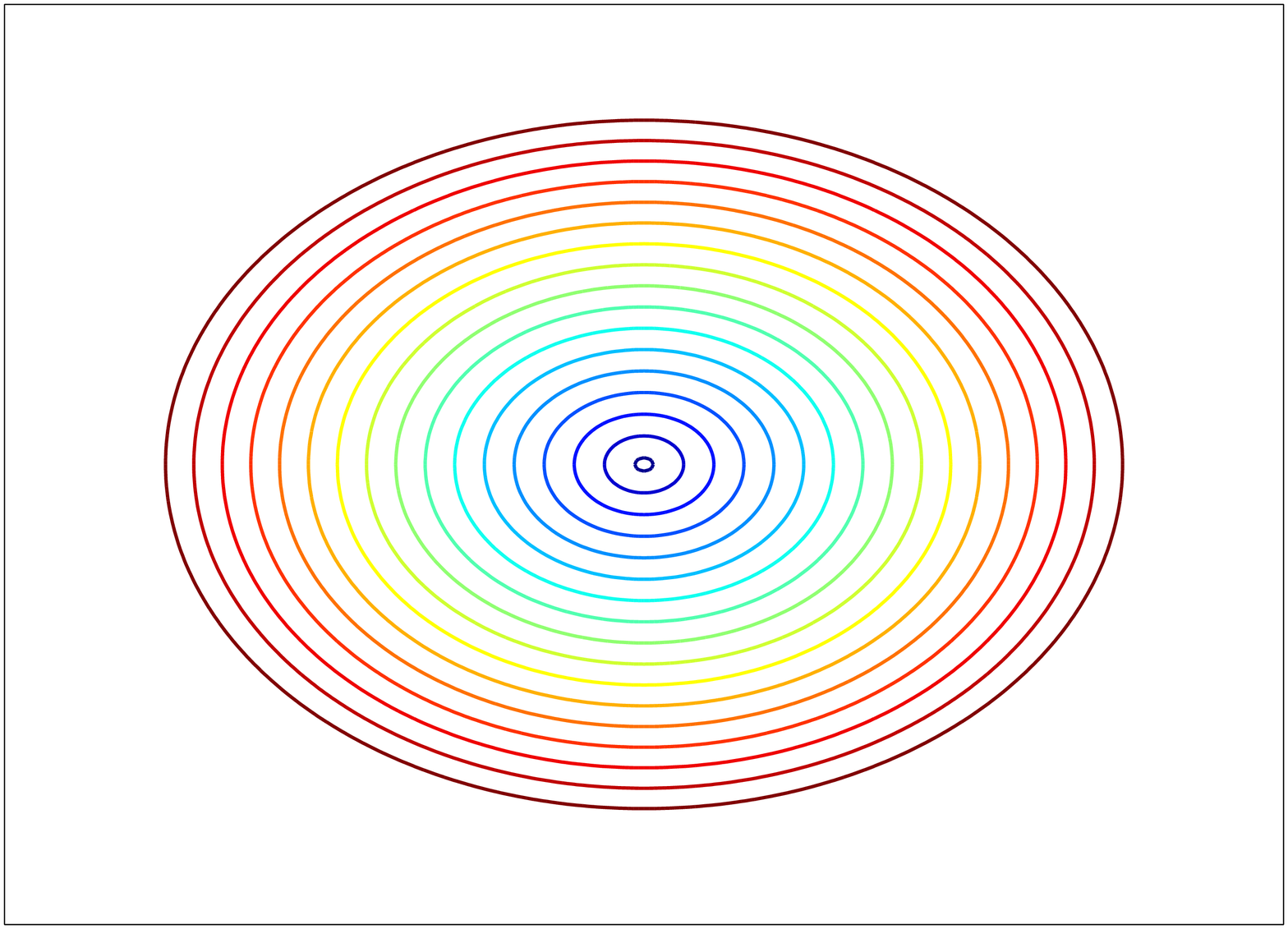}&
\includegraphics[width=0.3\textwidth]{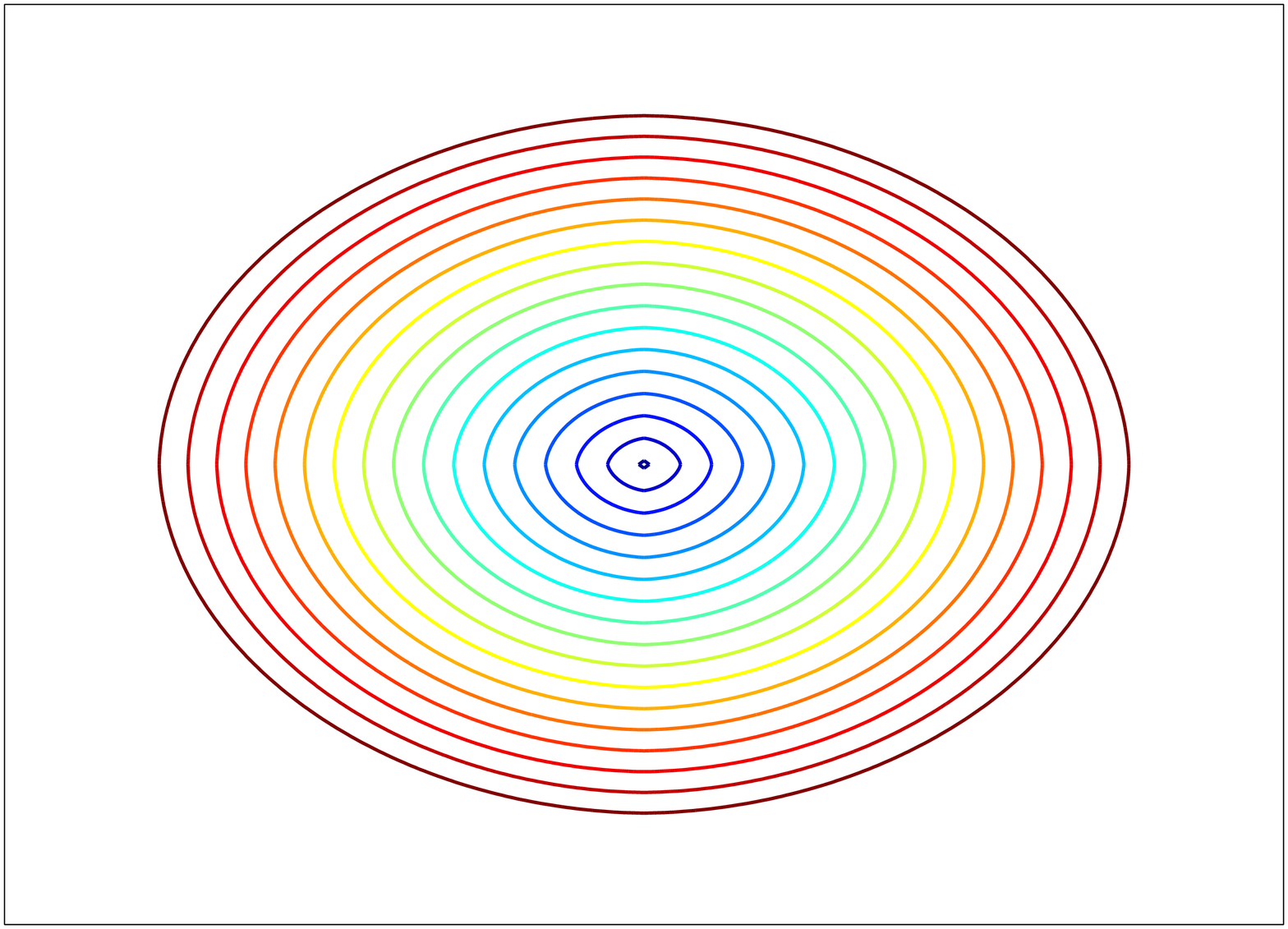}\tabularnewline
\end{tabular}
\par\end{centering}
\caption{Contour plots: (i) Left: True solution, (ii) Center: Schr\"odinger,
and (iii) Right: Fast sweeping}
\label{fig:pointsourceContourPlots} 
\end{figure}

The contour plots of the true solution and those obtained from Schr\"odinger
and fast sweeping are displayed below (figure ~\ref{fig:pointsourceContourPlots}).
We can immediately observe the similarity of our solution with the
true solution. We do observe smoother isocontours in our Schr\"odinger
method relative to fast sweeping.

\subsubsection{Comparison with the fast sweeping}
In order to verify the accuracy of our
technique, we compared our solution with fast sweeping for the following set of examples, using the
latter as the ground truth as the true solution is not available in closed-form.

\noindent \textbf{Example 2:} In this example we solved the eikonal equation from a point source
located at $(x_{0},y_{0})=(0,0)$ for the following forcing function
\begin{equation}
f(x,y)=1+2(e^{-2((x+0.05)^{2}+(y+0.05)^{2})}-e^{-2((x-0.05)^{2}+(y-0.05)^{2})})\end{equation}
 on a $2D$ grid consisting of points between $(-0.125,-0.125)$
and $(0.125,0.125)$ with a grid width of $\frac{1}{2^{10}}$. 
We ran our method for $6$ iterations with $\hbar$ set at
$0.015$ and fast sweeping for $15$ iterations sufficient for both
techniques to converge. When we calculated the percentage error for
the Schr\"odinger according to equation~\ref{eqn:percentError} (with
fast sweeping as the ground truth), the error was just around $1.245\%$.
The percentage error and maximum difference between the fast sweeping and Schr\"odinger solutions
after each iteration are adumbrated in Table~(\ref{table:eikonalexp2}).

\begin{figure}[ht!] 
\begin{centering}
\begin{tabular}{cc}
\includegraphics[width=0.4\textwidth]{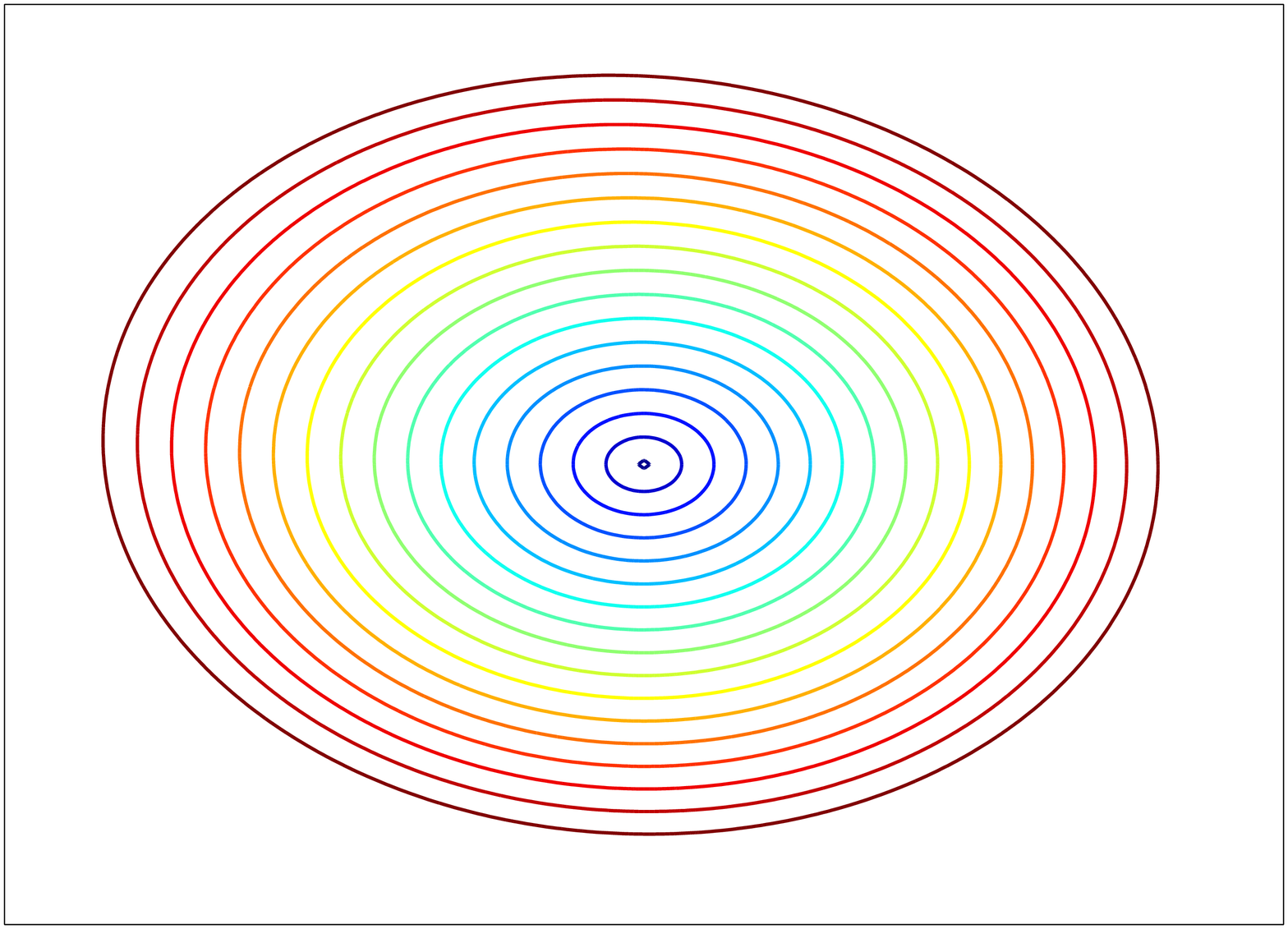}&
\includegraphics[width=0.4\textwidth]{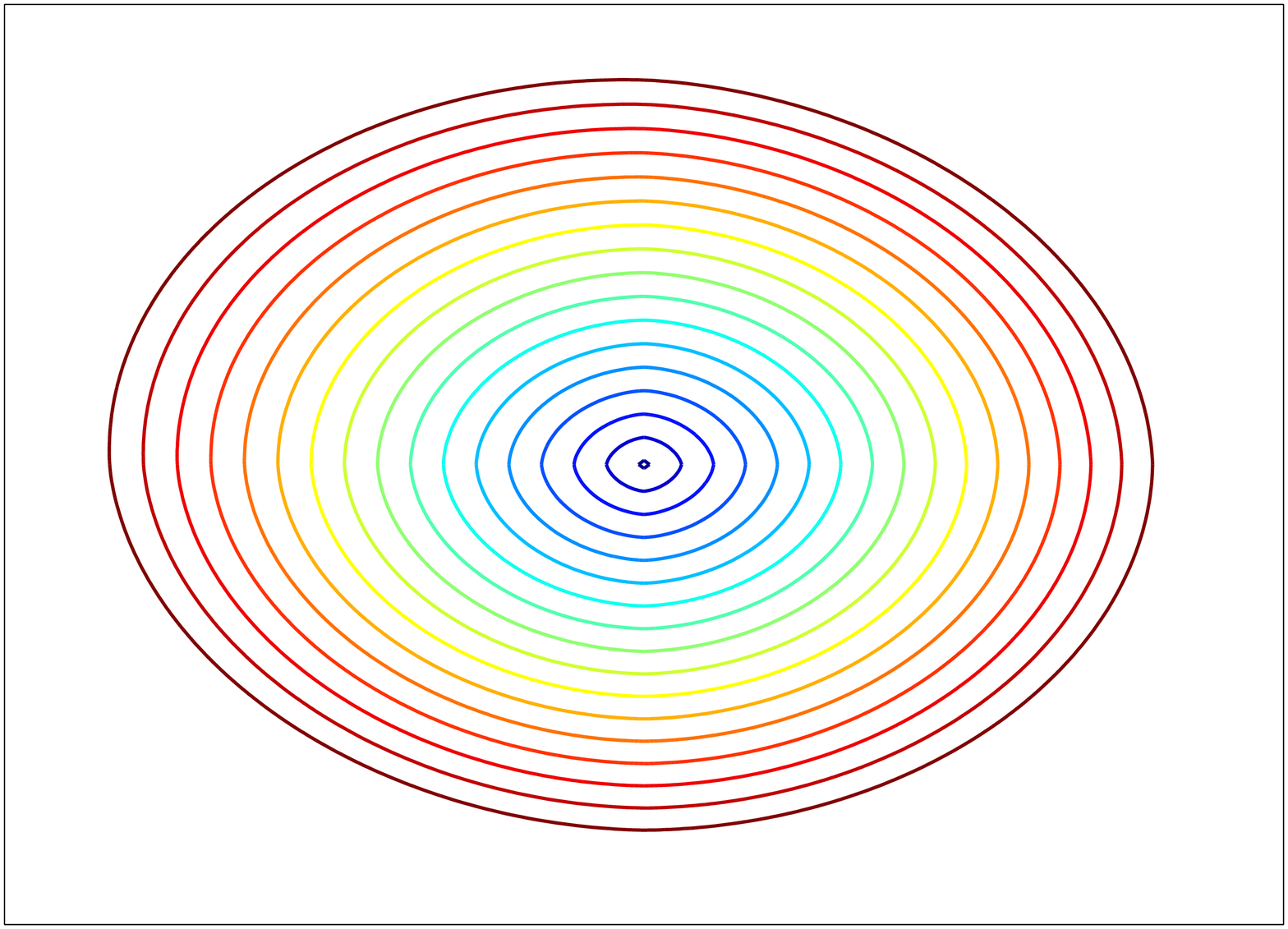}\tabularnewline 
\end{tabular}
\par\end{centering}
\caption{Contour plots: (i) Left: Schr\"odinger, (ii) Right: Fast sweeping}
\label{fig:pointsourceContourPlotsExp2} 
\end{figure}

\noindent \begin{center}
\begin{table}[ht!] 
\centering{}
\caption{Percentage error for the Schr\"odinger method in comparison to fast sweeping.}
\label{table:eikonalexp2}
\begin{tabular}{|c|c|c|}
\hline 
Iter  & $ \% error$  & max diff\tabularnewline
\hline 
1 & 1.144632 & 0.008694 \tabularnewline
\hline 
2 & 1.269028 & 0.008274 \tabularnewline
\hline 
3 & 1.223836 & 0.005799 \tabularnewline
\hline 
4 & 1.246392 & 0.006560 \tabularnewline
\hline 
5 & 1.244885 & 0.006365 \tabularnewline
\hline 
6 & 1.245999 & 0.006413 \tabularnewline
\hline 
\end{tabular}
\end{table}
\end{center}
We believe that the fluctuations both in the percentage error and the maximum difference are due
to repeated approximations of the integration involved in the convolution with discrete convolution
and summation, but nevertheless stabilized after 6 iterations. The contour plots shown in Figure~\ref{fig:pointsourceContourPlotsExp2}
clearly demonstrate the similarities between these methods.

\noindent \textbf{Example 3:} Here we solved the eikonal equation 
for the sinusoidal forcing function
\begin{equation}
f(x,y)=1+\sin(\pi(x-0.05))\sin(\pi(y+0.05))
\end{equation}
 on the same $2D$ grid as in the previous example. We randomly
chose $4$ grid locations namely,
\begin{center}
$\{0,0\},\{ 0.0488,0.0977\},\{-0.0244,-0.0732 \},\{ 0.0293,-0.0391\}$
\end{center}
as data locations and ran our method for $6$ iterations with $\hbar$ set at
$0.0085$ and ran fast sweeping for $15$ iterations. The percentage error between the
Schr\"odinger solution (after 6 iterations) and fast sweeping 
was $4.537 \%$ with the maximum absolute difference between them being $0.0109$.

\begin{figure}[ht!]
\begin{centering}
\begin{tabular}{cc}
\includegraphics[width=0.4\textwidth]{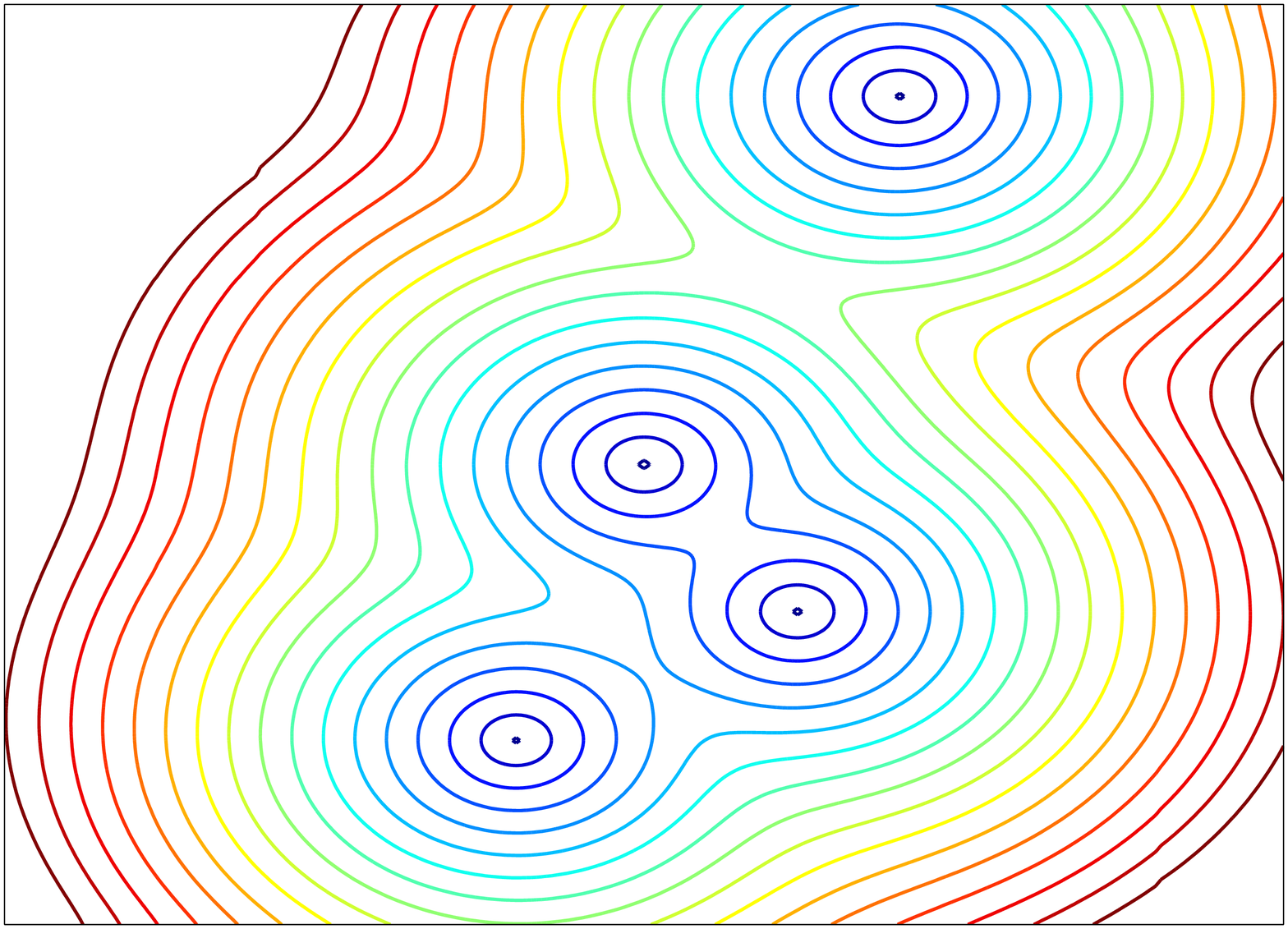}&
\includegraphics[width=0.4\textwidth]{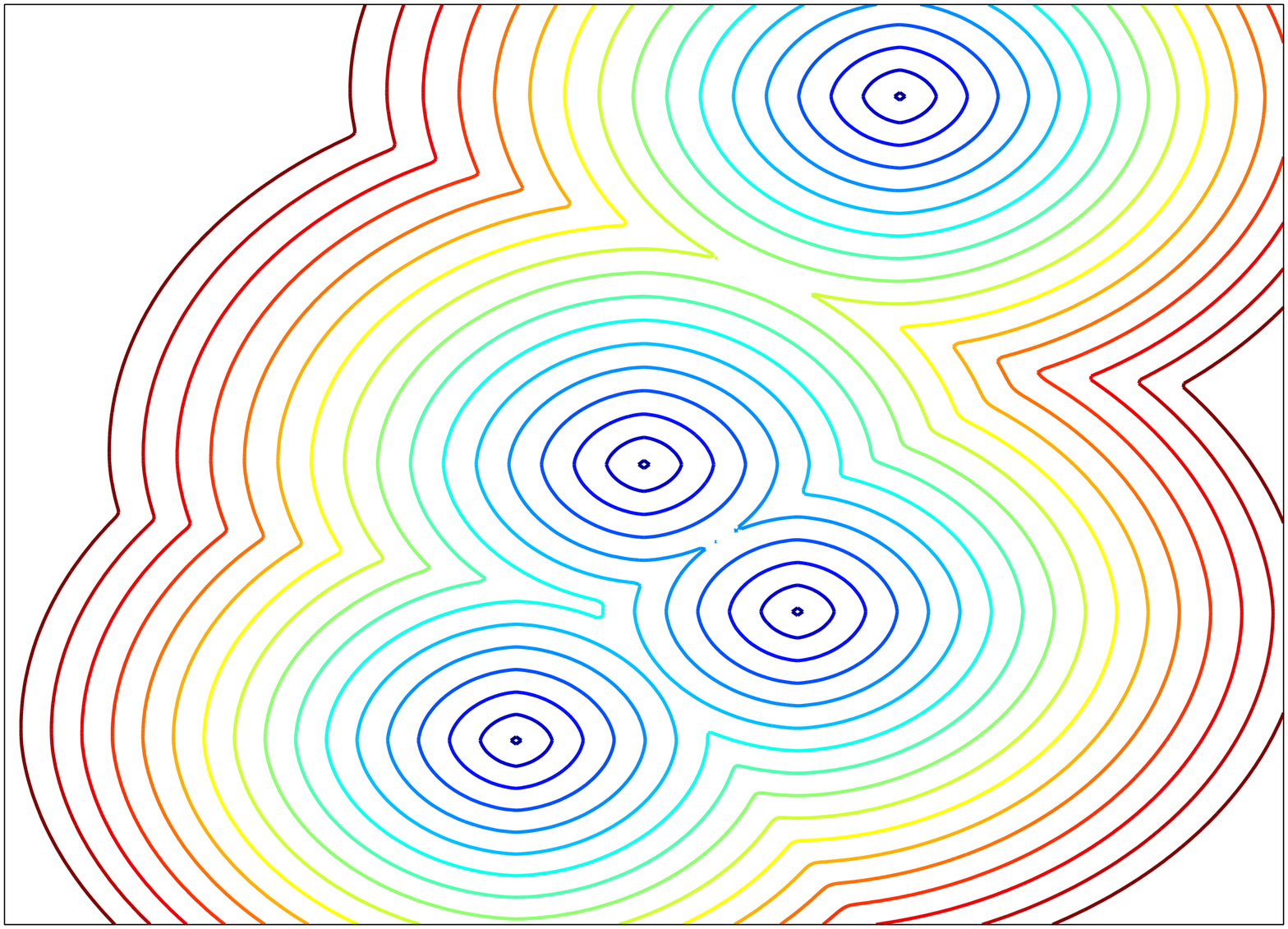} \tabularnewline
\end{tabular}
\par\end{centering}
\caption{Contour plots: (i) Left: Schr\"odinger, (ii) Right: Fast sweeping}
\label{fig:4pointsContourPlotsExp} 
\end{figure}

The contour plots are shown in Figure~(\ref{fig:4pointsContourPlotsExp}). Notice that the Schr\"odinger contours
are more smoother in comparison to the fast sweeping contours.

\noindent \textbf{Example 4:} Here we compared with fast sweeping
on a larger $2D$ grid consisting of points between $(-5,-5)$ and $(5,5)$
with a grid width of $0.25$. We again considered the sinusoidal forcing function
\begin{equation}
f(x,y)=1+0.3\sin(\pi(x+1))\sin(\pi(y-2))
\end{equation}
and chose 4 grid locations namely $\{0,0\}, \{1,1\},\{-2,-3\},\{3,-4\} $
as data locations. Notice that the Green's function $G$ and $\tilde{G}$
goes to zero exponentially faster for grid locations away from zero for small values of $\hbar$.
Hence for a grid location say $(-4,4)$ which is reasonably far away from 0, the value of the Green's function
say at $\hbar = 0.001$ may be zero even when we use a large number of precision bits $p$. This problem can be easily circumvented by first scaling down the entire grid by a factor $\tau$, computing the solution $S^{\ast}$ on the smaller denser grid and then rescaling it back again by $\tau$ to obtain the actual solution. It is worth emphasizing that scaling down the grid is tantamount to scaling down the forcing function as clearly seen from the fast sweeping method. In fast sweeping \cite{Zhao05}, the solution $S^{\ast}$ is computed using the quantity $f_{i,j}\delta$ where $f_{i,j}$ is the value of forcing function at the $(i,j)^{th}$ grid location and $\delta$ is the grid width. Hence scaling down $\delta$ by a factor of $\tau$ is equivalent to fixing $\delta$ and scaling down $f$ by $\tau$. Since the eikonal equation (\ref{eikonalEq}) is linear in $f$, computing the solution for a scaled down $f$--equivalent to a scaled down grid--and then rescaling it back again is guaranteed to give the actual solution.

\begin{figure}[ht!]
\begin{centering}
\begin{tabular}{cc}
\includegraphics[width=0.4\textwidth]{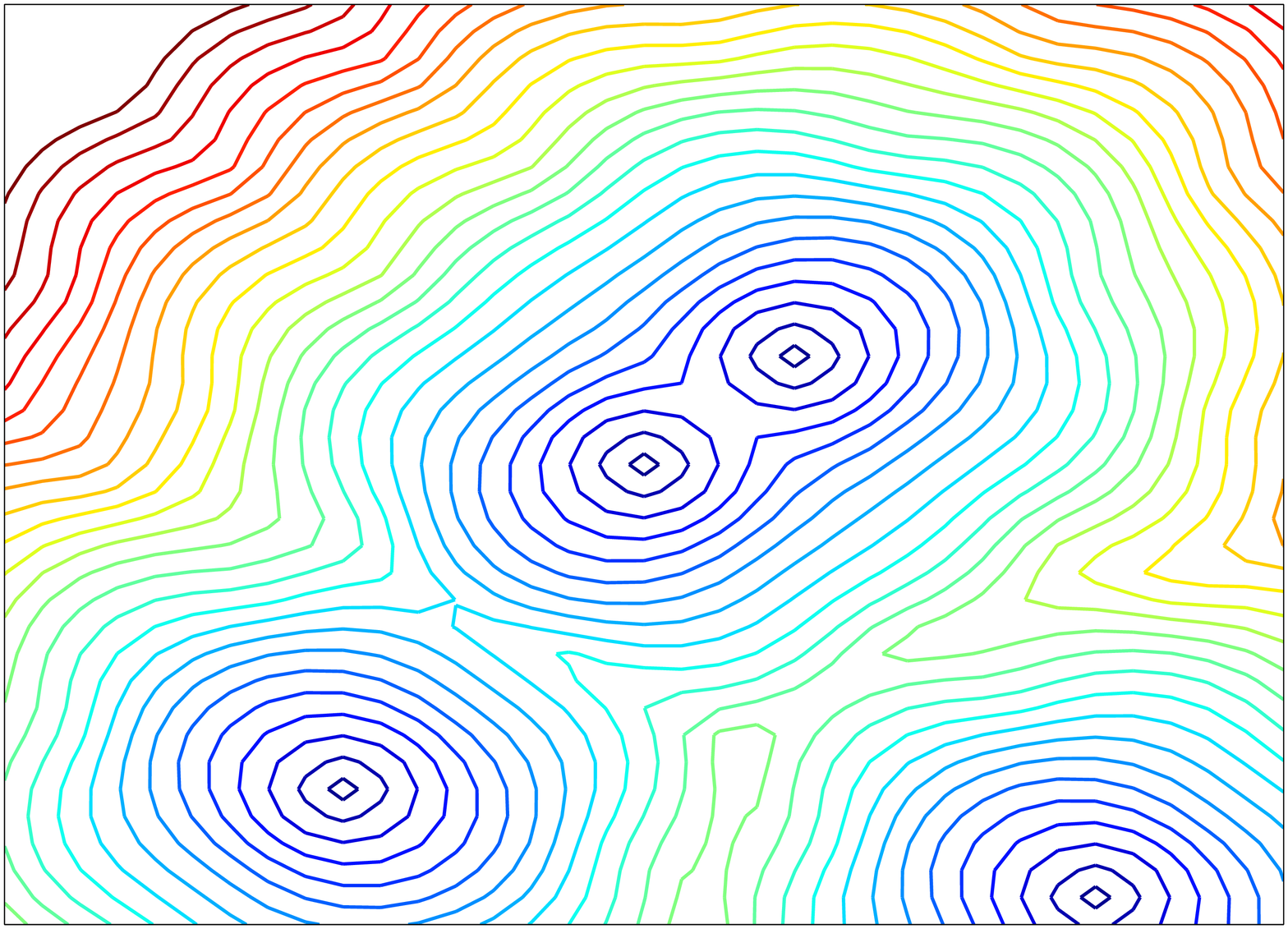}&
\includegraphics[width=0.4\textwidth]{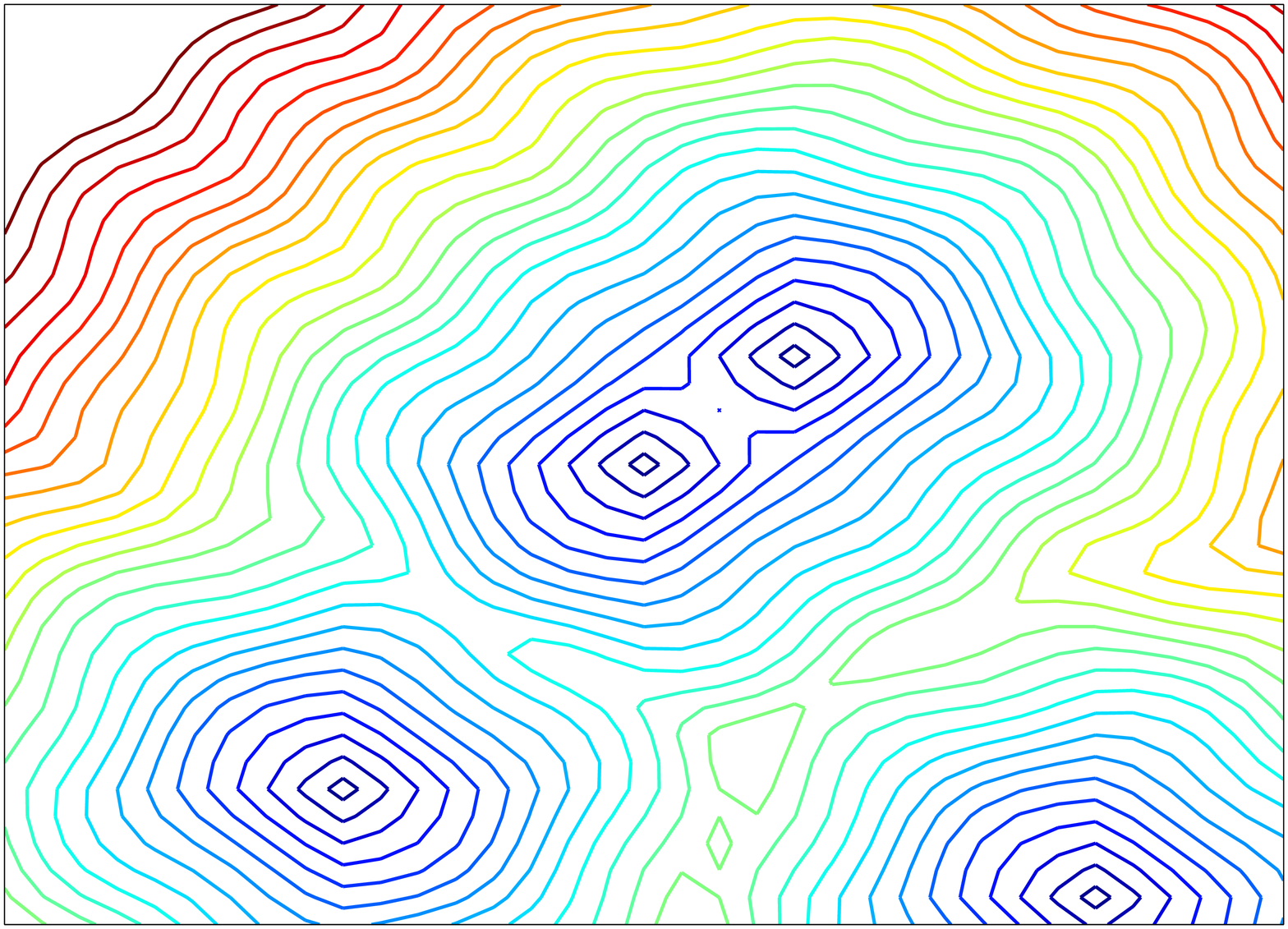}\tabularnewline 
\end{tabular}
\par\end{centering}
\caption{Contour plots: (i) Left: Schr\"odinger, (ii) Right: Fast sweeping}
\label{fig:4pointsContourPlotsLargerGrid} 
\end{figure}

The factor $\tau$ can be set to any desired quantity. For the current experiment we set $\tau = 100$, $\hbar = 0.001$ and ran our method for 6 iterations. Fast sweeping was run for 15 iterations. The percentage error between these methods was about $3.165 \%$. The contour plots are shown in Figure~\ref{fig:4pointsContourPlotsLargerGrid}. Again, the contours obtained from the Schr\"odinger are more smoother than those obtained from fast sweeping.

\section{Discretization approach to solve the linear system with applications to path planning and shape from shading}
\label{sec:discretizationApproach}
As seen above, the perturbation technique requires repeated convolution of the Green's function $G$ with the solution from the previous iteration $\phi_{i-1}$. As this convolution does not carry a closed form solution in general, we approximate the continuous convolution by its discrete counterpart computed via FFT. The error incurred from this discrete approximation tends to pile up with iteration. To circumvent this issue, a possible direct route to solving the screened Poisson equation in (\ref{eq:phiEquation}) is by discretization of the linear operator ($\nabla^2$) on a standard grid where the Laplacian is approximated by standard finite differences. Solving the discretized screened Poisson equation is far less complicated to implement than the fast marching or fast sweeping methods needed for directly solving (\ref{eikonalEq}). In addition, the computational complexity for implementing (\ref{eq:phiEquation}) can match these algorithms since $O(N)$ sparse, linear system solvers
are available \cite{Saad03}. This comes from the fact that a finite
difference approximation to (\ref{eq:phiEquation}), using a standard
five-point Laplacian stencil, simply results in a sparse linear system
of the form \cite{Demmel97}:
\begin{align}
&
\underbrace{
\left[\begin{array}{cccc}
L_{N}+2I_{N} & -I_{N} & \cdots & 0\\
-I_{N} & \ddots & \ddots & \vdots\\
\vdots & \ddots & \ddots & -I_{N}\\
0 &  & -I_{N} & L_{N}+2I_{N}
\end{array}\right]
+
\begin{bmatrix}
f(X_1)&0&\cdots&0\\
0&f(X_2)&\cdots&0\\
\vdots&\vdots&\ddots&\vdots\\
0&0&\cdots&f(X_N)\\
\end{bmatrix}
}_{A}
\underbrace{\left[\begin{array}{c}
\phi_{1}\\
\phi_{2}\\
\vdots\\
\phi_{N}
\end{array}\right]}_{x}=
\underbrace{\left[\begin{array}{c}
\delta_{kron}(\x_1)\\
\vdots\\
\delta_{kron}(\x_0)\\
\vdots\\
\delta_{kron}(\x_N)
\end{array}\right]}_{b} \nonumber \\
&
\label{eq:sparse-system}
\end{align}
where $N$ represents the number of grid points and $L_{N}$ is a tri-diagonal
block of the form
\[
L_{N}=\left[\begin{array}{cccc}
2 & -1 &  & 0\\
-1 & \ddots & \ddots\\
 & \ddots & \ddots & -1\\
0 &  & -1 & 2
\end{array}\right].
\]
The function $\delta_{kron}(\x)$ in (\ref{eq:sparse-system}) is defined in (\ref{eq:deltakron}) which is supported and takes the value one
only at the point-set locations $\{\y_k\}_{k=1}^K$. Though better complexities are
achievable using multigrid \cite{Saad03} solvers, one can just as
well address many problems in a satisfactory manner using direct,
sparse solvers, like MATLAB's $A\backslash b$---an approach we adopt
for the experiments in the present paper. 

\subsection{Path planning} 
In the pioneering contribution of \cite{Kimmel01}, Kimmel and Sethian  developed one of the earliest applications of the eikonal equation
to path planning. By finding a solution $S^{\ast}(\x)$ (referred to as the
value function in the path planning context) to the eikonal equation
in (\ref{eikonalEq}), we immediately recover the minimum
cost to go from a source location $\x_{0}$ in the state space to any
other point $\x$ (in the state space). Here, we impose the boundary
condition $S^{\ast}(\x_0)=0$, and consider $f(\x)$ as the cost to travel
through location $\x$ (higher the value, the more costly) and prescribe
it as a strictly positive function. The function $f(\x)$
function is set to high values for undesirable travel regions
for the optimal path, and very small values for favorable travel areas
(and set to one on the source point). In comparison to other popular path planning techniques like potential
field methods \cite{Khatib86}, the value function is an example of
a navigation function---a potential field free of local minima. 

Given a scalar field solution, the optimal paths are determined by
gradient descent on $S^{\ast}(\x)$, and is typically referred to as backtracking.
The backtracking procedure can be formulated as an ordinary differential
equation 
\begin{equation}
\dot{\x}=-\frac{\nabla S^{\ast}(\x(t))}{\|\nabla S^{\ast}(\x(t))\|},
\label{eq:PathPlanning-backTrackOnS}
\end{equation}
whose solution $\x(t)$ is the reconstructed path from a fixed target
location $\hat{\x}$. We typically terminate the gradient backtracking
procedure at some $t$ value such that $\left\Vert \x(t)-\x_0\right\Vert $ $<\epsilon$,
i.e. we get arbitrarily close the source point, for some small $\epsilon>0$.
By construction, the backtracking on $S^{\ast}(\x)$ cannot get stuck in local
minima---an obvious proof by contradiction validates this claim if
one considers $S^{\ast}$ to be differentiable and have local minima $\nabla S^{\ast}=0$
at some $\x$, but $f(\x)>0$, and contradicts the eikonal equation
$\|\nabla S^{\ast}\|=f(\x)$. Although, in theory, $S^{\ast}(\x)$ can contain saddle
points, but usually this is not an issue in practice. 

\subsection{Anecdotal verification of path planning on complex mazes and extensions to vessel segmentation}
\label{sec:Experimental-Results}
\begin{figure}[t]
\begin{centering}
\begin{tabular}{cc}
\includegraphics[width=1.5in,height=1.25in]{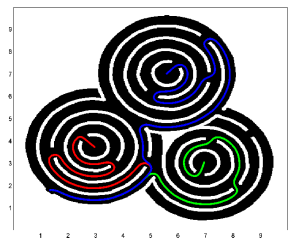}\tabularnewline 
\end{tabular}
\par\end{centering}
\caption{\label{fig:triCricMaze}
Optimal paths from various locations to common source. (Note: constraint
areas are in white and traversable regions in black.)}
\end{figure}

\begin{figure}[t]
\begin{centering}
\begin{tabular}{ccc}
\includegraphics[width=0.956in,height=1.047in]{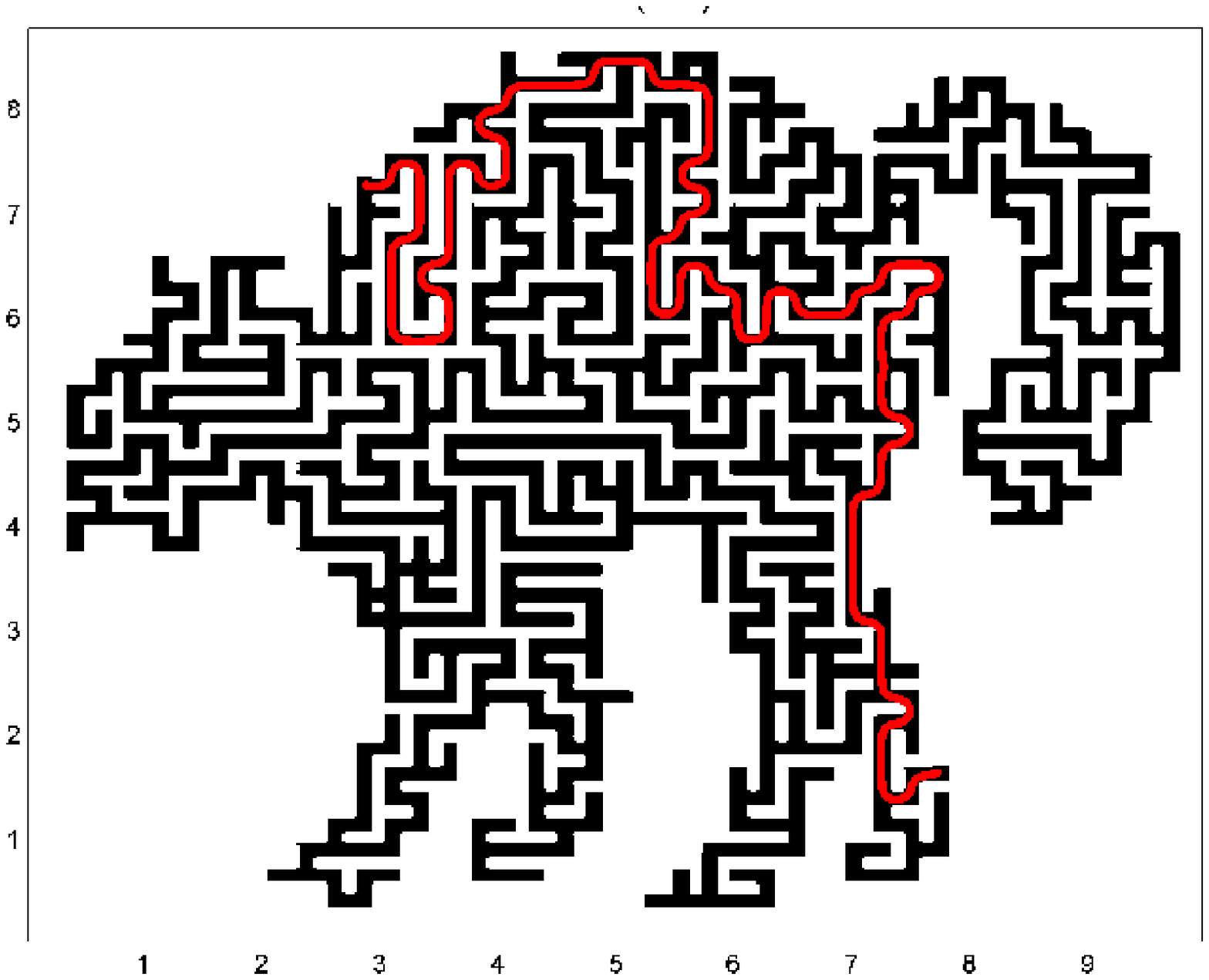}
& \includegraphics[width=0.956in,height=1.047in]{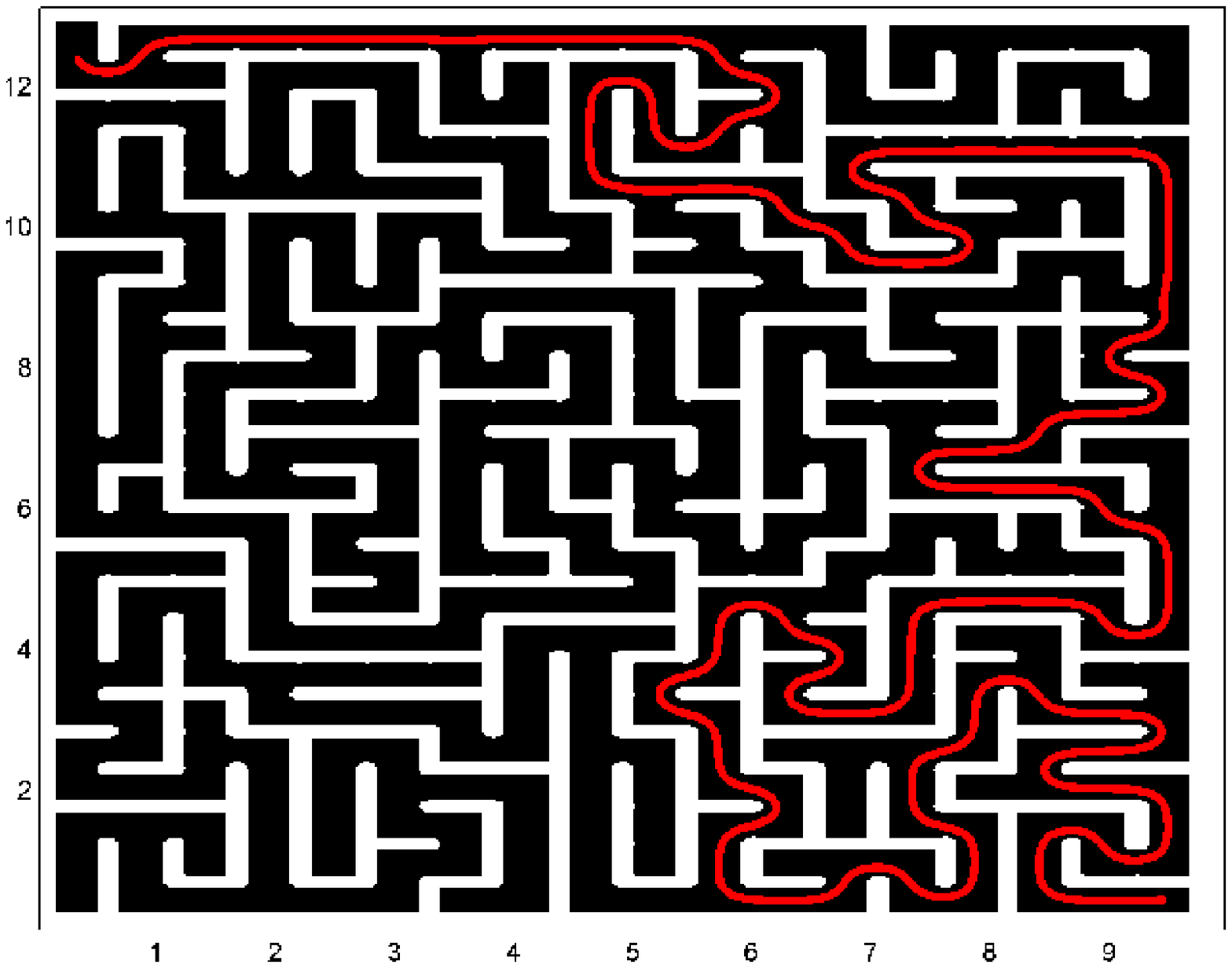}
&
\includegraphics[width=0.956in,height=1.047in]{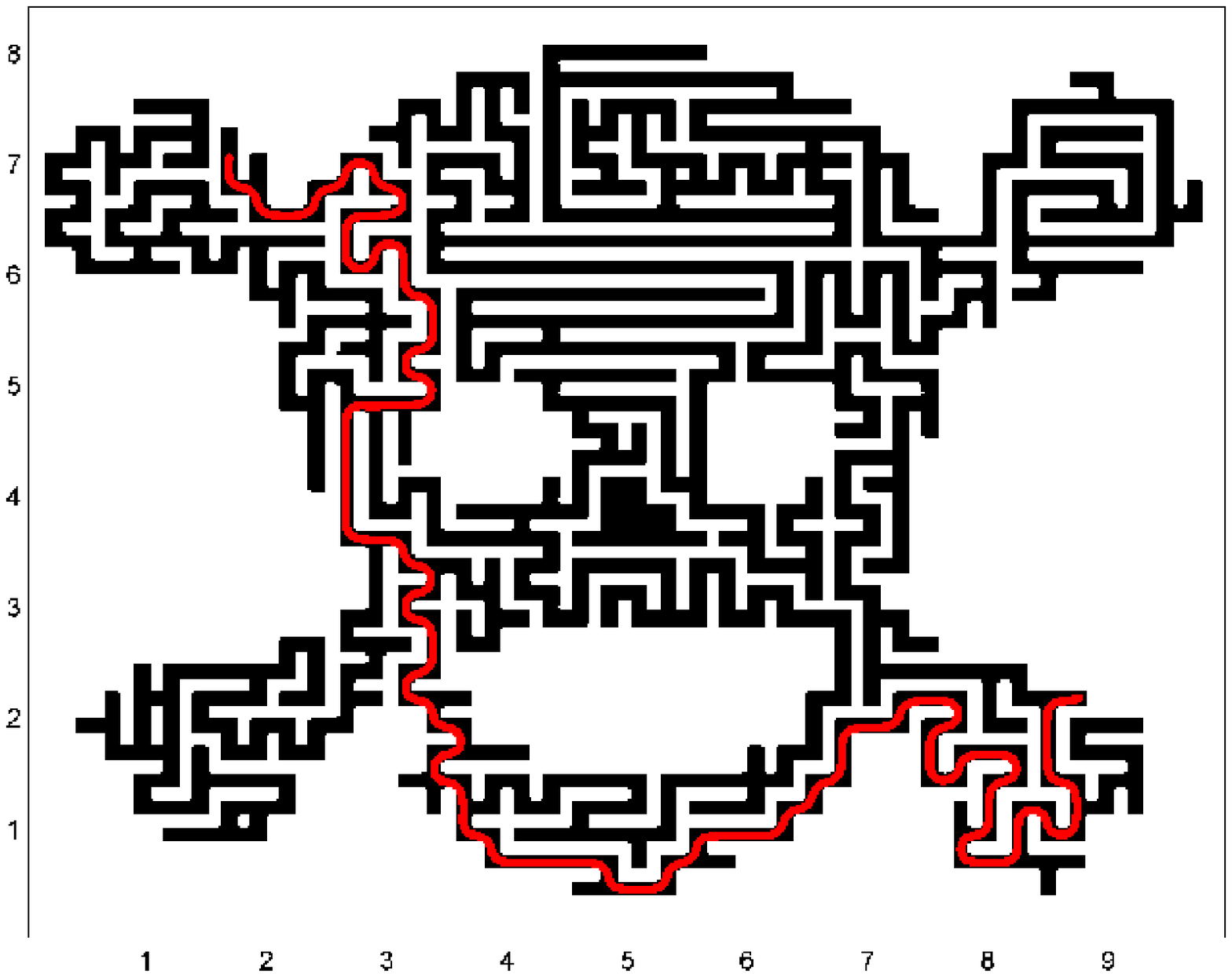}\tabularnewline 
(a) & (b) & (c)\tabularnewline
\end{tabular}
\par\end{centering}
\caption{\label{fig:catAndskull}(a) Schr\"{o}dinger's cat maze, (b) standard
maze, (c) skull maze. Multiple solutions are possible but only the
shortest path is chosen to be optimal. (Note: constraint areas are
in white and traversable regions in black.) }
\end{figure}

We applied our path planning approach to a variety of complex maze
images. We explicitly chose the maze grid sizes to be much larger
than the norm for recent publications that apply the eikonal equation
for path planning; these typical sizes are usually smaller than $100\times100$.
An objective juxtaposition of contemporary fast sweeping and marching
techniques, which require special discretization schemes, data structures,
sweep orders, etc., versus our approach presented here, clearly illustrates
the efficiency and simplicity of the later. Our framework reduces
path planning (a.k.a. all-pairs, shortest path or geodesic processing)
implementation to four straightforward steps:
\begin{enumerate}
\item Define $f(\x)$, which assigns a high cost to untraversable areas in
the grid and low cost to traversable locations. In the experiments
here, we simply let appropriately scaled versions of the maze images
be $f(\x)$, with \emph{white pixels representing boundaries} and \emph{black
pixels the possible solution paths}.
\item Select a source point on the grid and a small value for $\hbar$.
The solution to (\ref{eq:phiEquation}) will simultaneously recover
all shortest paths back to this source from any non-constrained region
in the grid.
\item Use standard finite differencing techniques to evaluate (\ref{eq:phiEquation}). 
This leads to a sparse, block tri-diagonal system which can be solved by a multitude
of linear system numerical packages. We simply use MATLAB's '\textbackslash{}'
operator. Recover approximate solution to eikonal by letting $S^{\ast}(\x)=-\hbar \log\phi(\x)$.
\item Backtrack to find the shortest path from any allowable grid location
to the source, i.e. use eq. (\ref{eq:PathPlanning-backTrackOnS})
to perform standard backtracking on $S^{\ast}(\x)$ .
\end{enumerate}
The grid sizes and execution times for several mazes are provided
in Table~\ref{tab:solveTime}. 
\begin{table}[t]
\begin{centering}
\caption{\label{tab:solveTime}Maze grids sizes and time to solve sparse system
for path planning. Our approach simply uses MATLAB's '\textbackslash{}'
operator to solve the shortest path problem, avoiding complex nuances
of discretization schemes and specialized data structures required
for fast marching.}
\begin{tabular}{|c|c|}
\hline 
2D Grid Dims. (No. of Points)  & $A\backslash b$ (sec.)\tabularnewline
\hline 
\hline 
{\footnotesize{$450\times450$ $\left(N=202,500\right)$, Fig. \ref{fig:triCricMaze}}} & {\footnotesize{0.79}}\tabularnewline
\hline 
{\footnotesize{$434\times493$ $\left(N=213,962\right)$, Fig. \ref{fig:catAndskull}(a)}} & {\footnotesize{1.02}}\tabularnewline
\hline 
{\footnotesize{$621\times473$ $\left(N=293,733\right)$, Fig. \ref{fig:catAndskull}(b)}} & {\footnotesize{1.22}}\tabularnewline
\hline 
{\footnotesize{$419\times496$ $\left(N=207,824\right)$, Fig. \ref{fig:catAndskull}(c)}} & {\footnotesize{0.77}}\tabularnewline
\hline 
\end{tabular}
\par\end{centering}
\end{table}
Notice that even for larger grids, our time to solve for $S^{\ast}(\x)$ is
on the order of a few seconds, and this is simply using the basic
sparse solver in MATLAB. The time complexity of MATLAB's direct solver
is $O\left(N^{1.5}\right)$, which makes our approach here slightly
slower than the optimal runtime. $O\left(N\right)$ is achievable
for sparse systems, such as ours, using multigrid methods, but we
have opted to showcase the simplicity of our implementation versus
pure speed. 
 Figure~\ref{fig:triCricMaze} illustrates the resulting
optimal paths from three different locations.
Figure~\ref{fig:catAndskull} illustrates our path planing
approach on a variety of mazes: (a) demonstrates path planning while
paying homage to Schr\"{o}dinger's cat, (b) is a traditional maze,
and (c) is a whimsical result on a skull maze. Notice in all these
mazes there are multiple solution paths back to the source, but only
the shortest path is chosen. 

\begin{figure}[t]
\begin{centering}
\begin{tabular}{ccc}
\includegraphics[width=0.945in,height=1.035in]{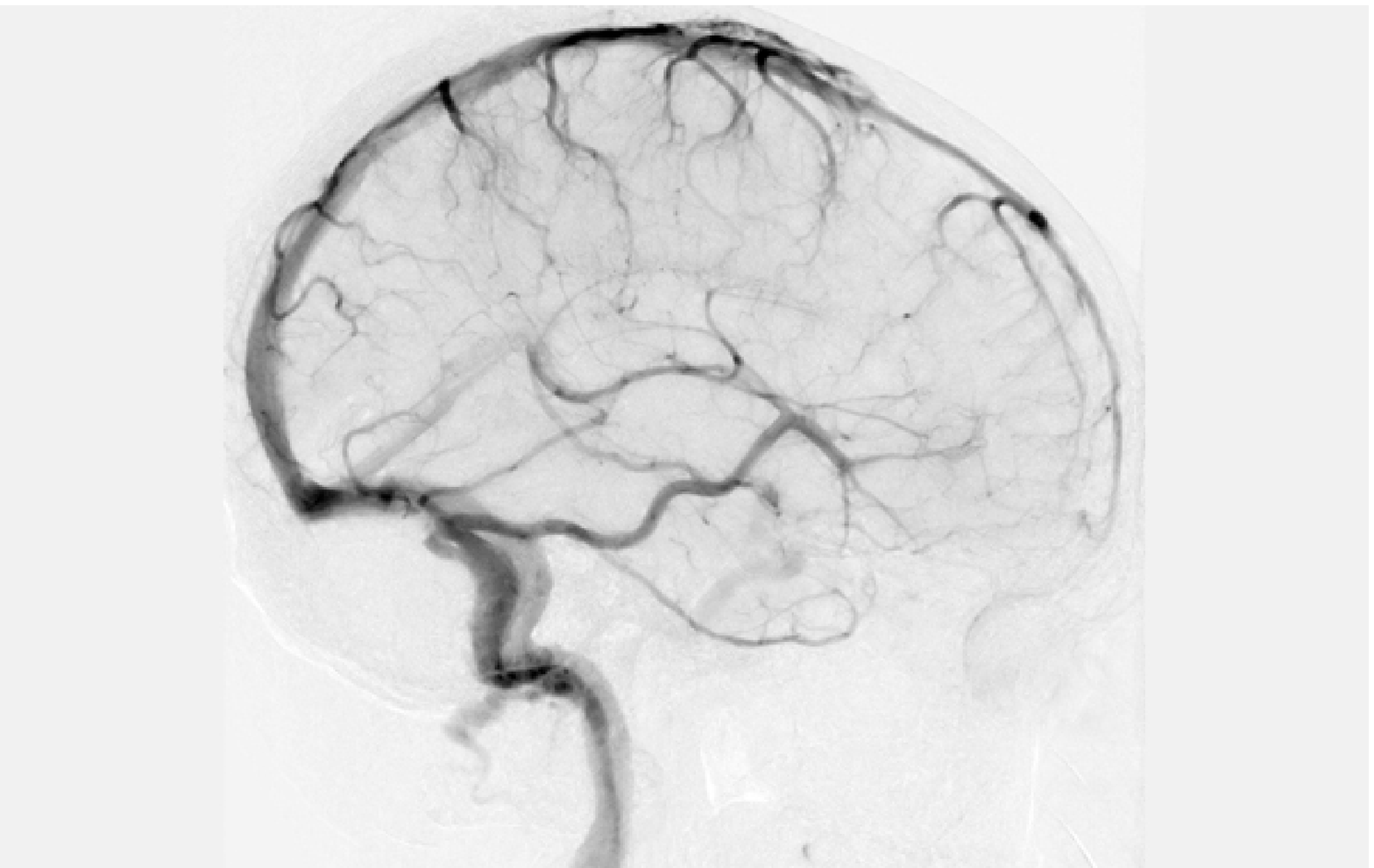} &
\includegraphics[width=0.945in,height=1.035in]{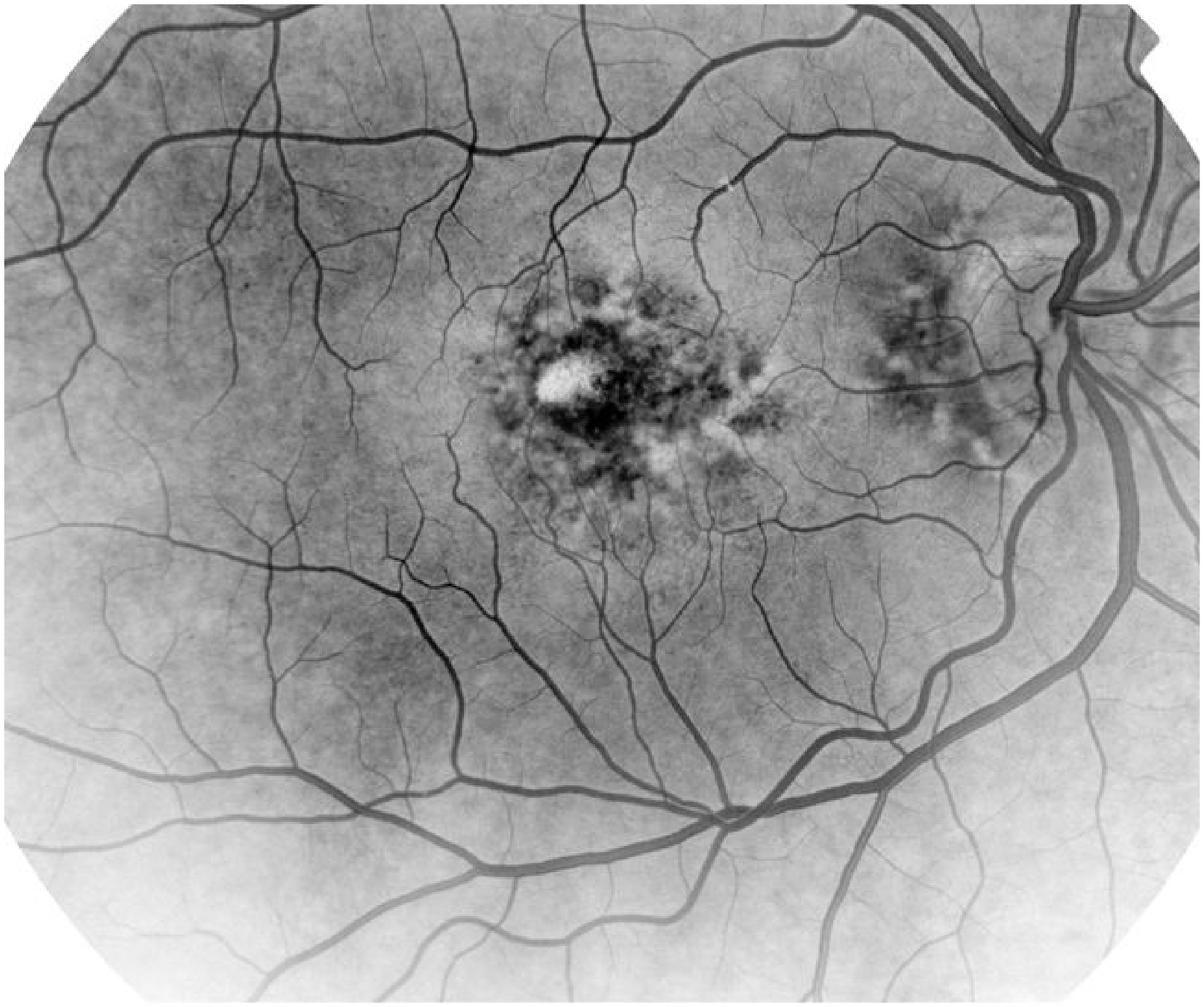} &
\includegraphics[width=0.945in,height=1.035in]{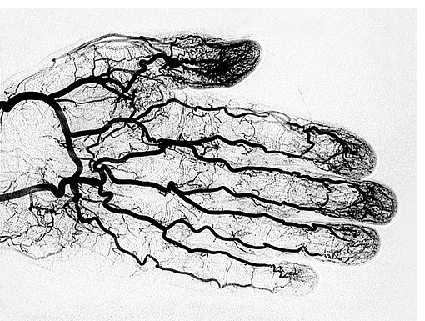}\tabularnewline 
\includegraphics[width=0.945in,height=1.035in]{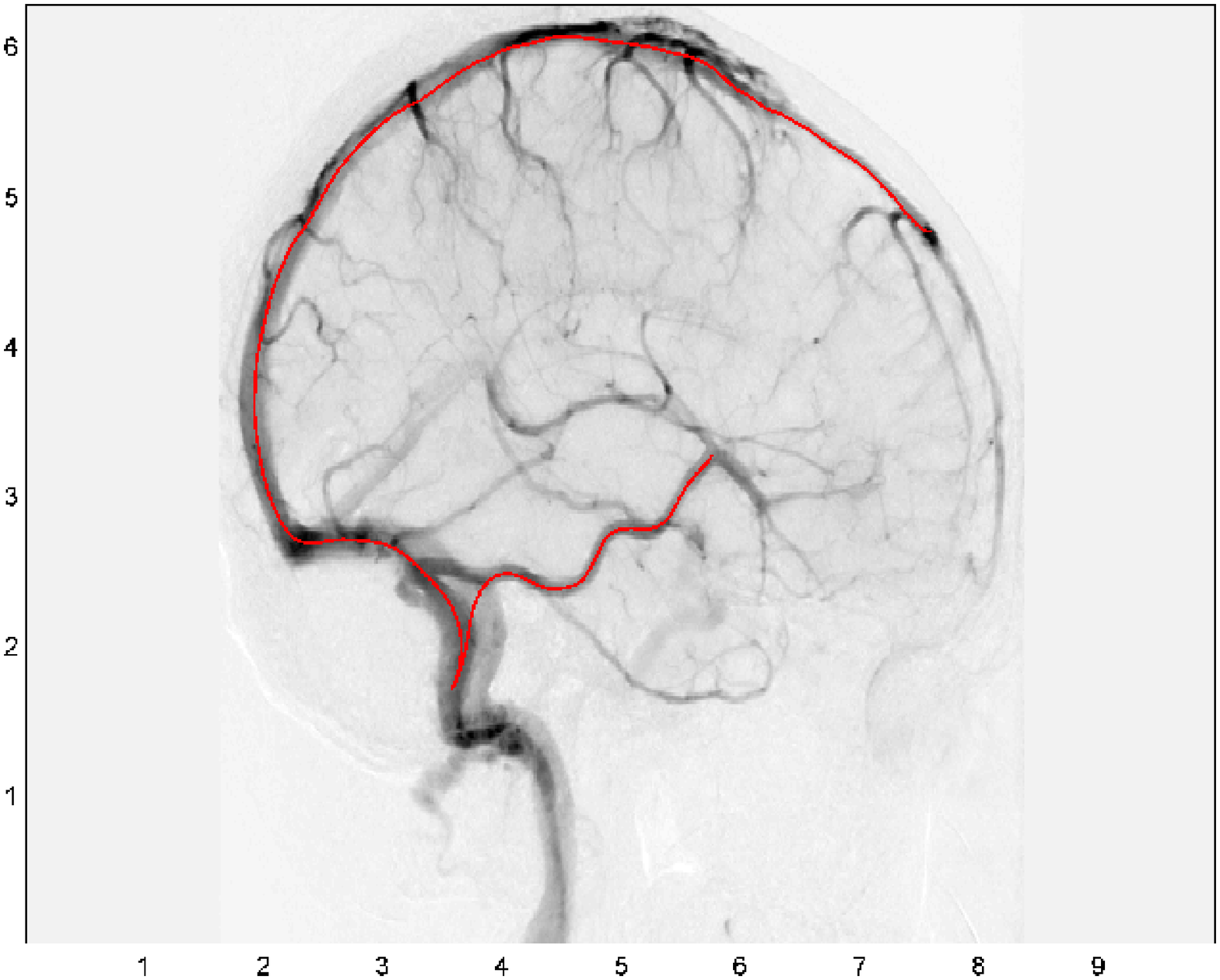} &
\includegraphics[width=0.945in,height=1.035in]{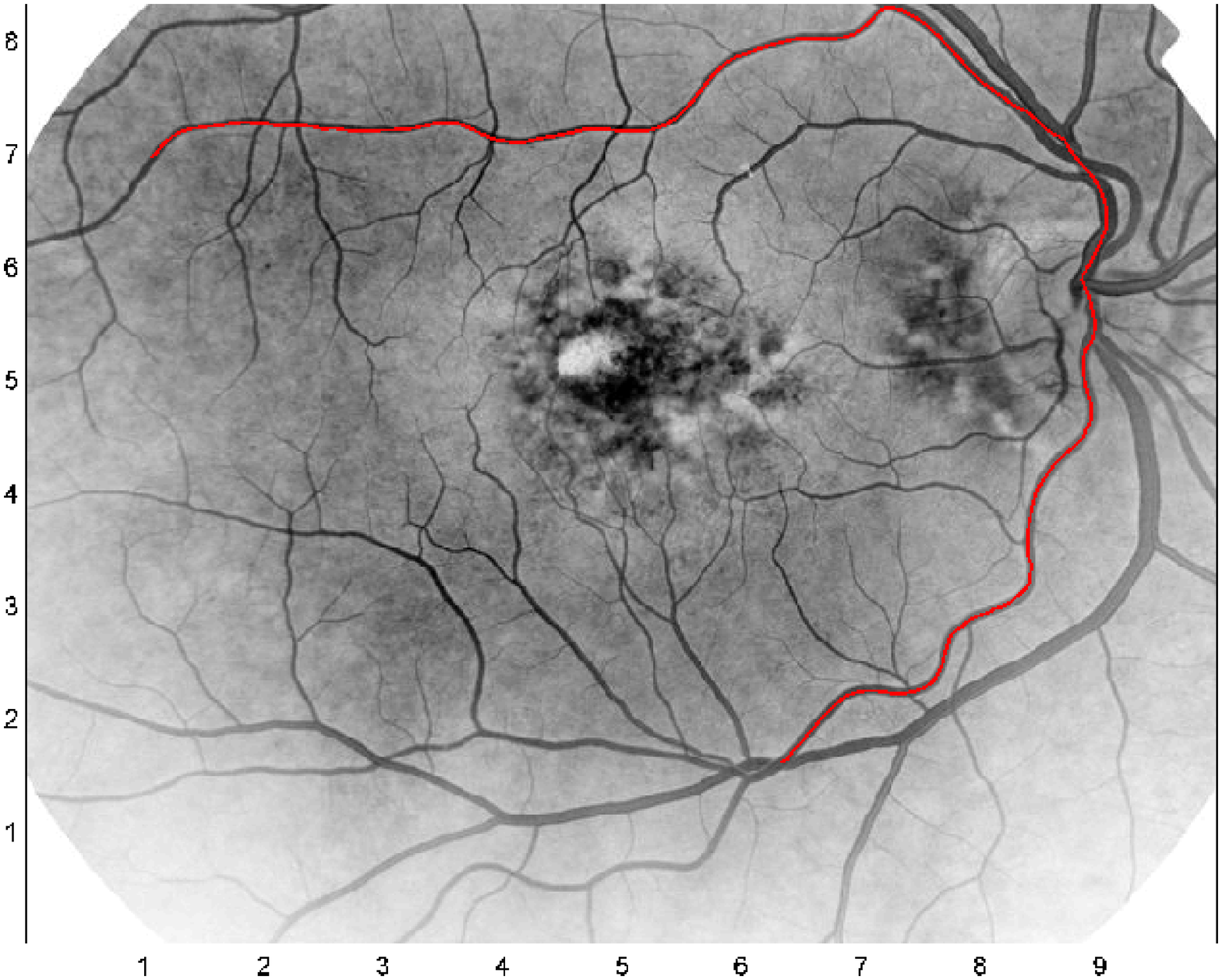} &
\includegraphics[width=0.945in,height=1.035in]{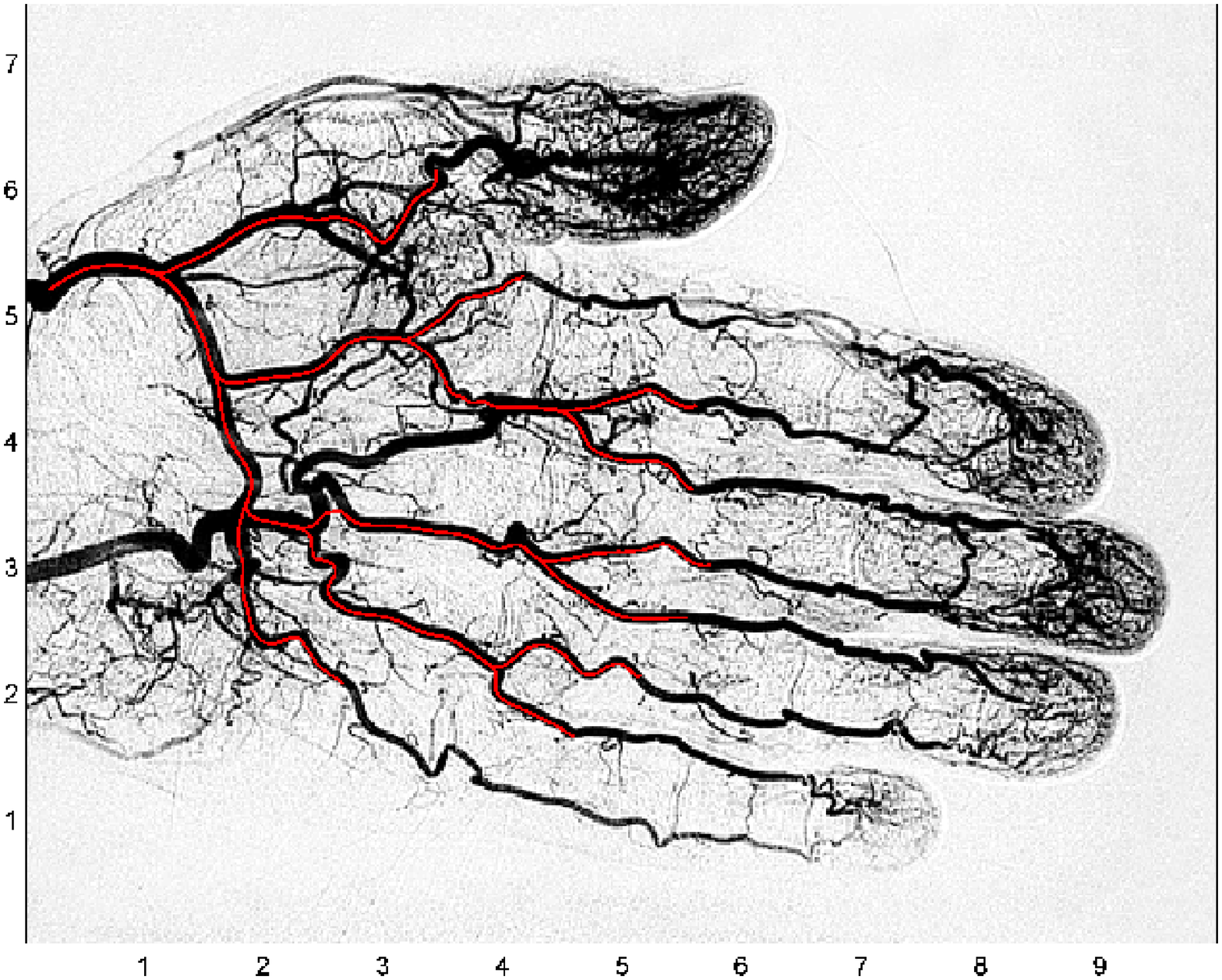}\tabularnewline 
\end{tabular}
\par\end{centering}
\caption{\label{fig:schVsFMandFS_vessSeg}Medical image vessel centerline extraction,
top row original images and bottom row our linear solution approach.
Paths under our linear systems approach are smooth due to inherent
viscosity-like behavior. They also produce segmentations where the
extracted vessels are centered on the blood vessels. Fast marching approach requires
additional constraints to achieve centerline extraction \cite{Cohen2001}. }
\end{figure}

The above discussed application of path planning can be readily extended
to centerline extraction from medical imagery of blood vessels. One
can view the image $I(\x)$ as a ``maze'' where we only want to travel
on the vessels in the image. Figure~\ref{fig:schVsFMandFS_vessSeg},
column (a) illustrates three example medical images: eye, brain, and
hand. Columns (b) showcases our results, while (c) provides comparative
analysis against fast sweeping. Notice that our solution naturally
generates smoother centerline segmentations, which is a natural consequence
of having a built-in, viscosity-like term in (\ref{Sequation}).
Whereas, the fast marching and fast sweeping methods tend to have
sharper transitions in the paths and deviate from the center---viscosity
solutions can be used to alleviate this, but are not organic to the
formulation like ours. In fact, it has been shown that additional
constraints have to be incorporated to ensure fast marching approaches
extract the centerline \cite{Cohen2001}. Again, we stress the simplicity
and ease-of-use of this approach, with only one free parameter $\hbar$,
making it a viable option for many path planning related applications,
such as robotic navigation, optimal manipulation, and vessel extraction
in medical images. 

\subsection{Shape-from-Shading}
\label{sec:Apps-Shape-from-Shading}
Shape-from-shading has long been a popular problem domain for computer
vision, having the primary objective of recovering the scalar height
field from a single image. Solution approaches utilizing the eikonal
equation have been known since the early 80's \cite{Bruss82}, and
have continually improved upon through the advent of fast sweeping
and fast marching methods \cite{Kimmel01,Pardos06}. The standard
forward image model assuming a Lambertian reflectance model generates
the luminance via inner product of the surface normal, $\n(\x)$, with
the light source direction, $\d$., i.e. $P(\x)=\langle \n(\x),\d\rangle$.
For example, if we assume a vertical lighting direction $\d=[0,\,0,\,1]^{T}$,
we get the imaging operator 
\[
P(\x)=\frac{1}{\sqrt{\left\Vert \nabla S^{\ast}(\x)\right\Vert ^{2}+1}},
\]
where $S^{\ast}(\x)$ is the desired scalar height field we wish to recover.
This can be obtained by solving the standard eikonal equation (\ref{eikonalEq}) with 
\begin{equation}
f(\x)=\sqrt{\frac{1}{P(\x)^{2}}-1}
\label{eq:SfS-forcingFunc}
\end{equation}
and boundary conditions $S^{\ast}(\x_{i})=h_{i}$, i.e. we seed the boundary
conditions with the known heights $h_{i}$ at select grid locations
$X_{i}$. 

As we have detailed in previous sections, our formalism allows one
to address any general (non-linear) eikonal equation by solving the
linear screened Poisson equation in (\ref{eq:phiEquation}). One simply
needs to create the forcing function in (\ref{eq:SfS-forcingFunc}) and then solve the discretized sparse system as in (\ref{eq:sparse-system}).
This immediately yields the recovered height field. 

\subsection{Surface reconstruction via shape from shading}
\begin{figure}[t]
\begin{centering}
\begin{tabular}{ccc}
\includegraphics[width=0.882in,height=0.966in]{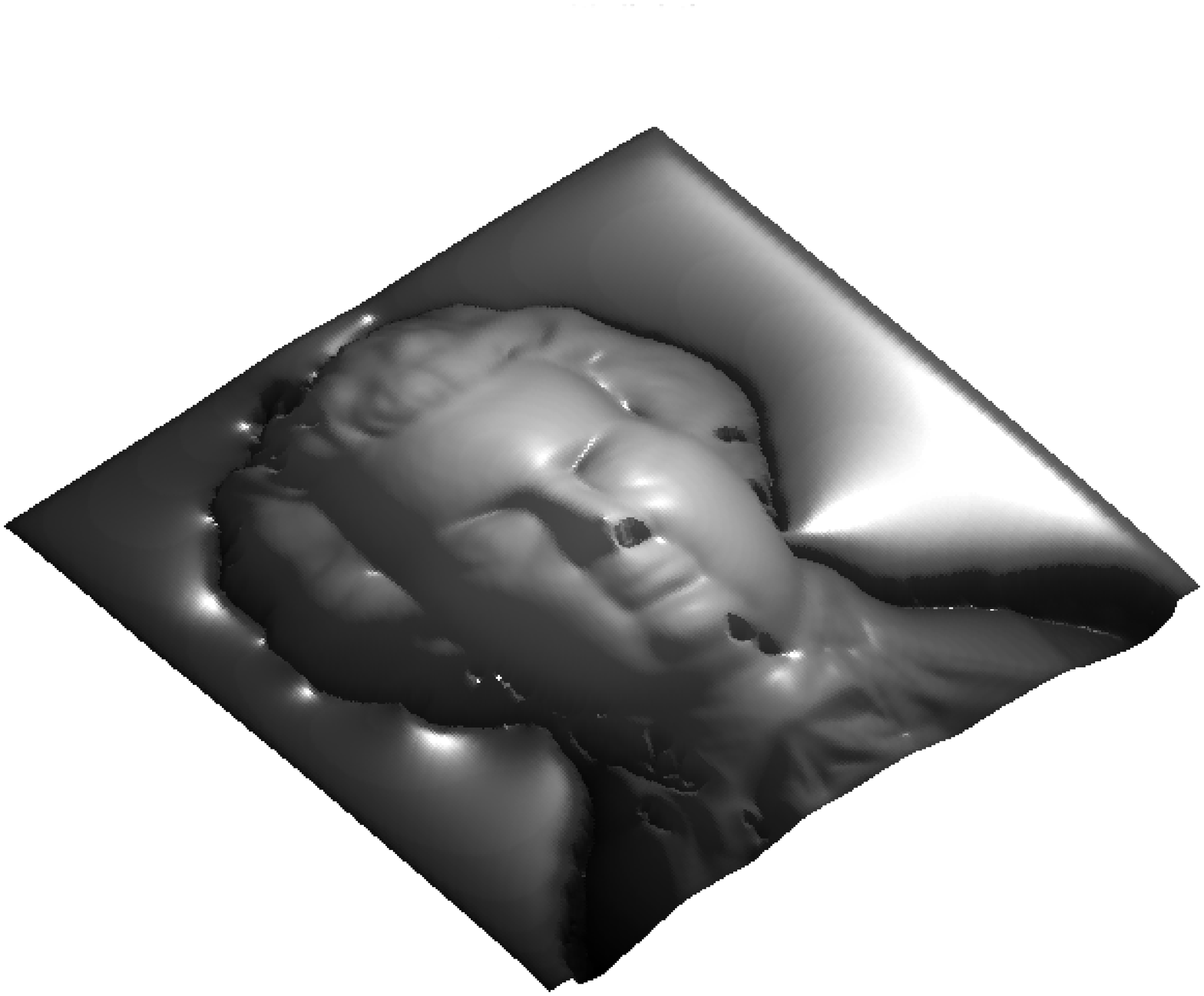} &
\includegraphics[width=0.882in,height=0.966in]{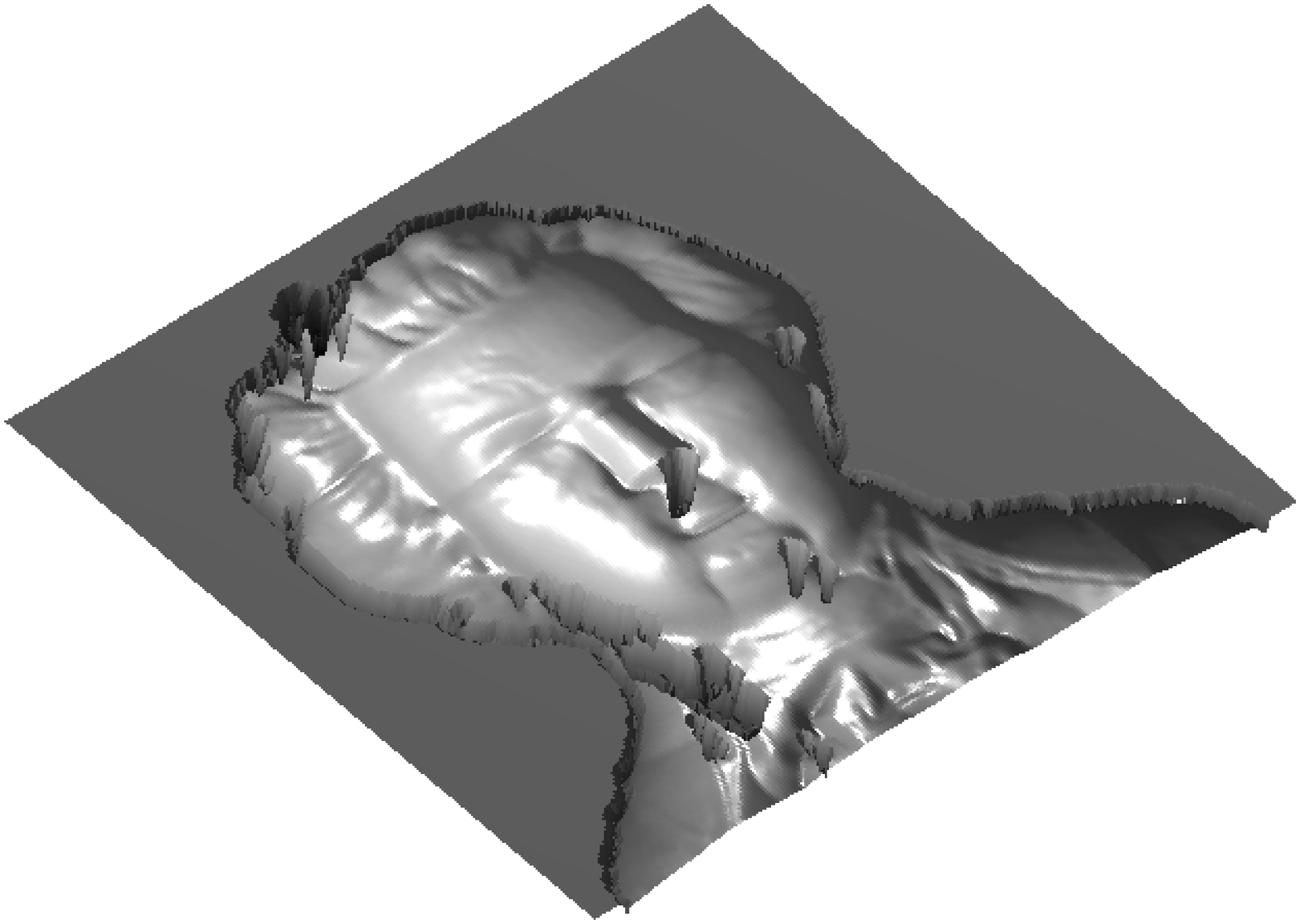}
&
\includegraphics[width=0.882in,height=0.966in]{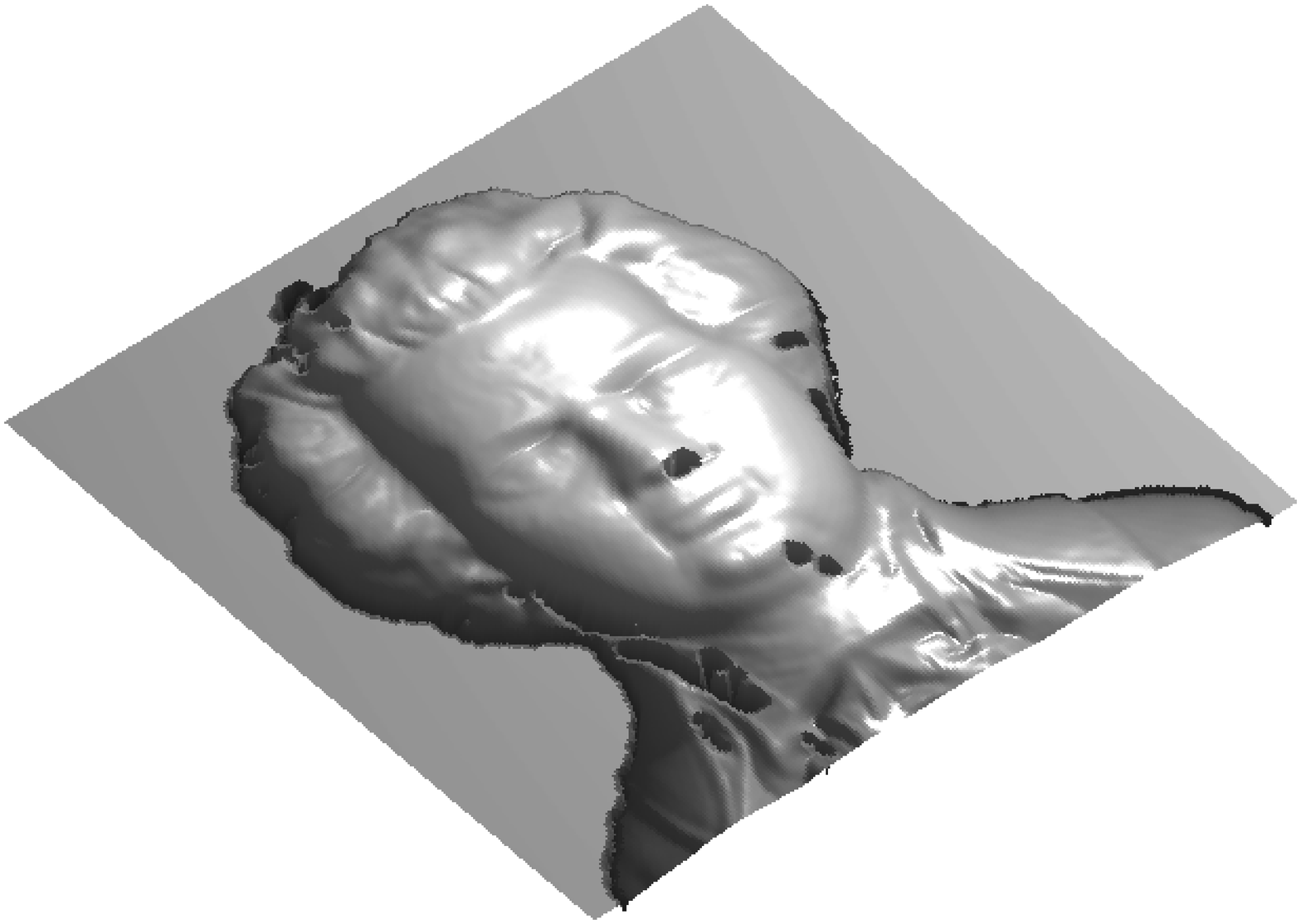}\tabularnewline 
Error: 0.524438 & Error: 0.713825 & Error: 0.654674\tabularnewline
\includegraphics[width=0.882in,height=0.966in]{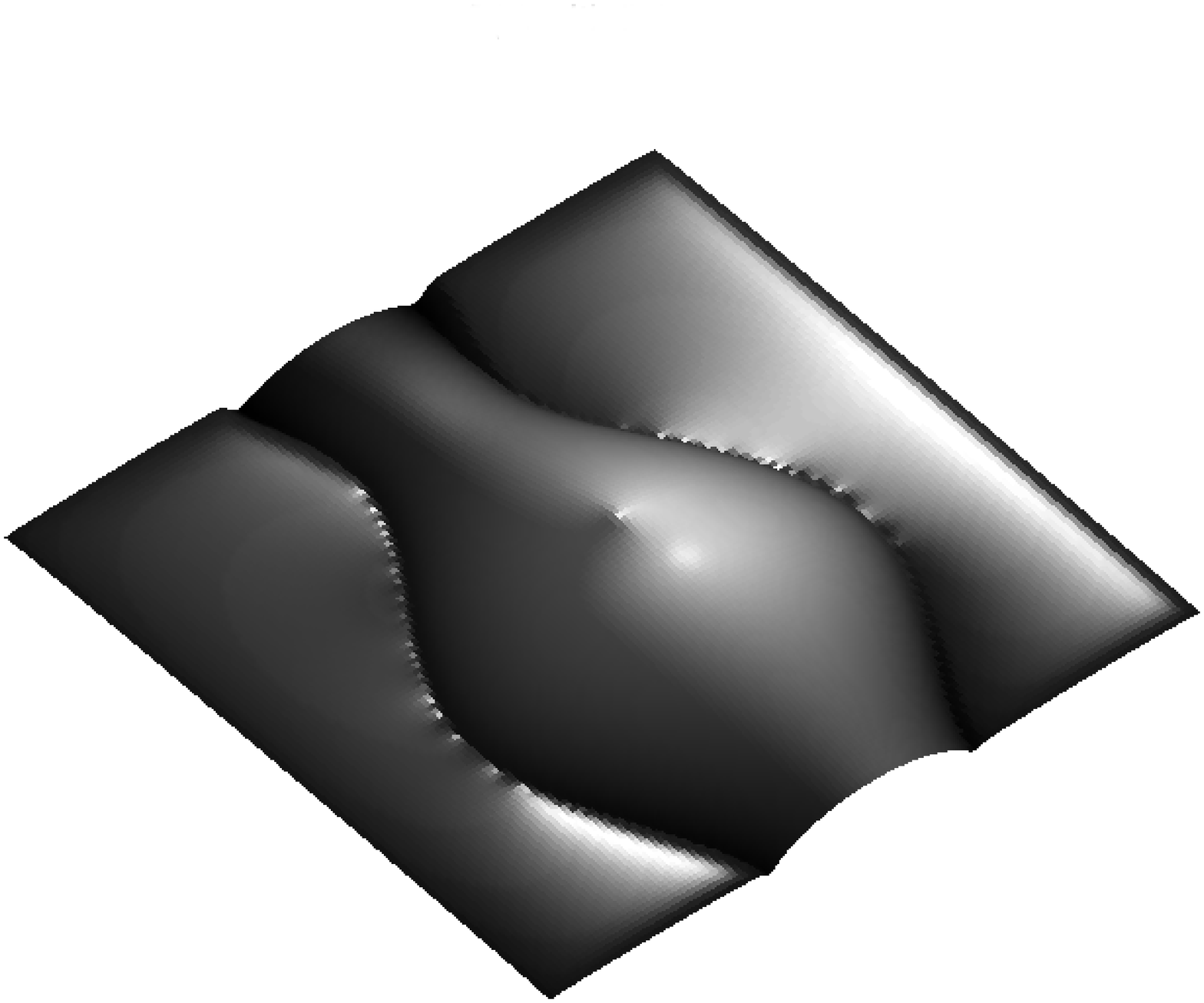} &
\includegraphics[width=0.882in,height=0.966in]{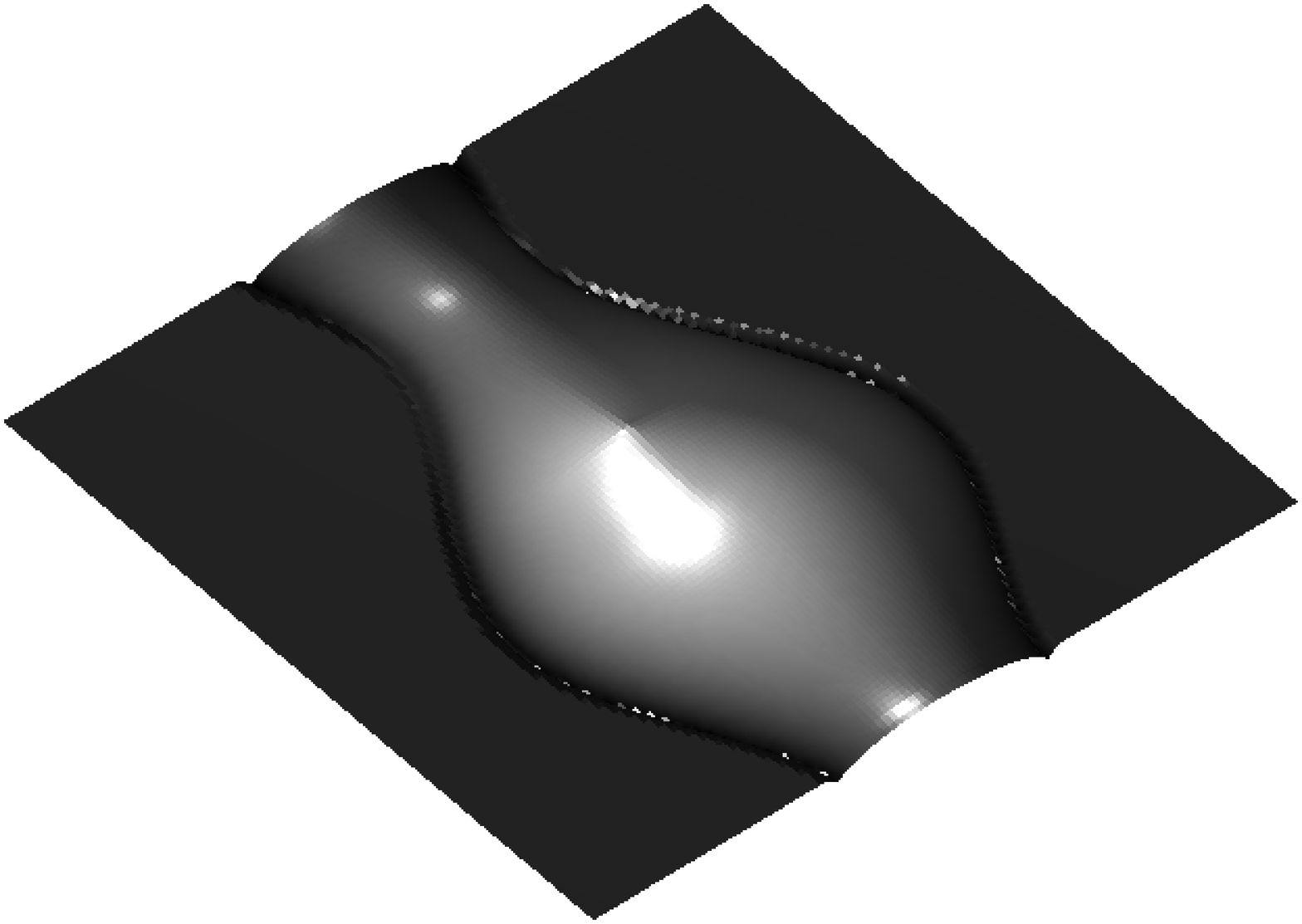} &
\includegraphics[width=0.882in,height=0.966in]{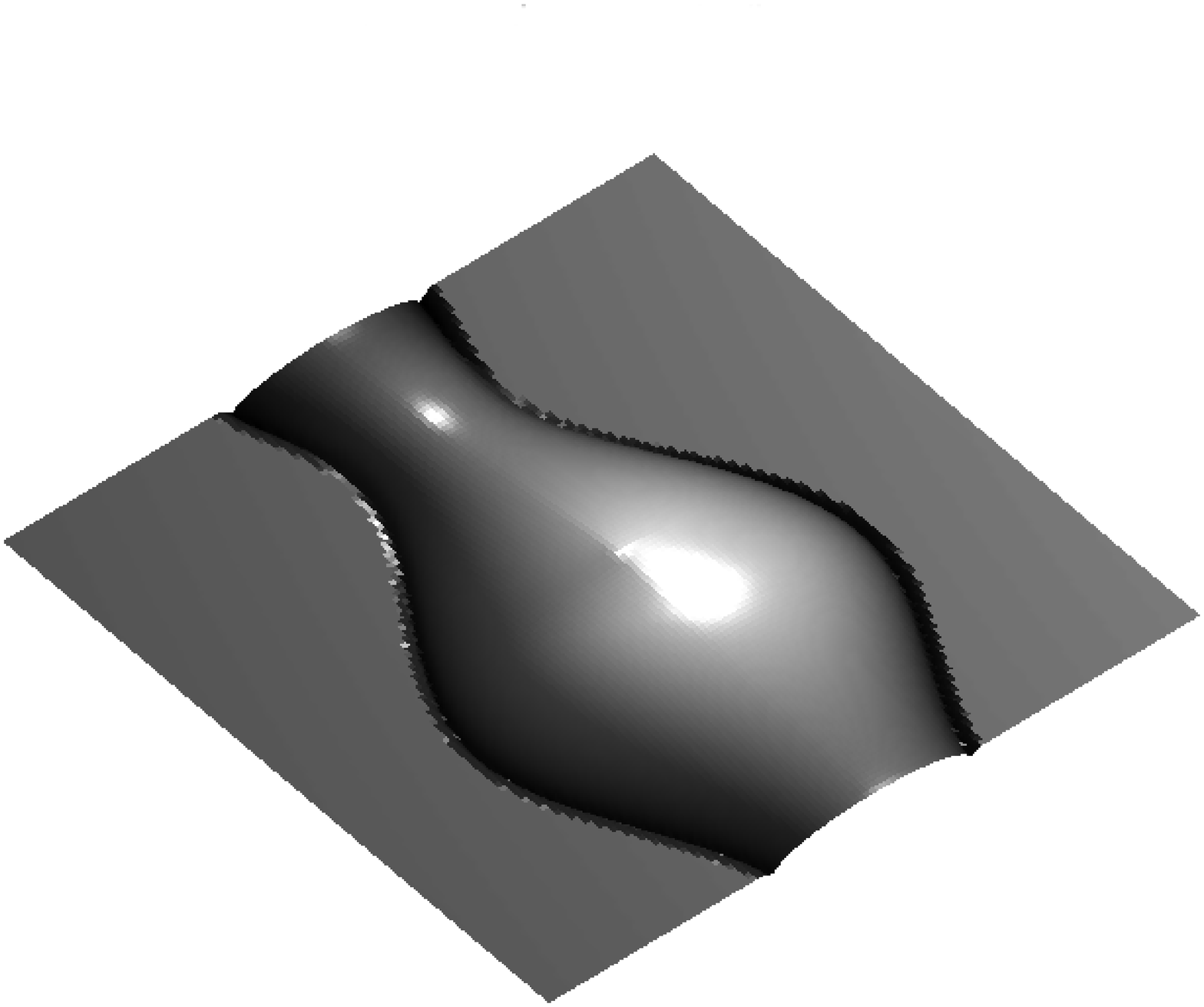}\tabularnewline 
Error: 0.203321 & Error: 0.237820 & Error: 0.199234\tabularnewline
(a) & (b) & (c)\tabularnewline
\end{tabular}
\par\end{centering}
\caption{\label{fig:schVsFMandFS_SfS}Shape-from-shading surface reconstruction,
per column: (a) our linear solution approach, (b) fast marching, (c)
fast sweeping. Based on the gradient magnitude error (from the true
surface), our approach linear systems approach is better or at least
highly competitive. }
\end{figure}

For shape from shading, we validated height recovery on two common
images that often used in the literature: Mozart and a vase. Figure~\ref{fig:schVsFMandFS_SfS}
illustrates the recovered surfaces using our method, (a), fast marching,
(b), and fast sweeping, (c). Under each image we also list the error
of the reconstruction from the known ground truth height field. The
error was computed by comparing the true mean gradient magnitudes
versus those estimated from the recovered $S^{\ast}(\x)$. 

The validation shows that our method is competitive with both fast
marching and fast sweeping, all the while retaining the efficiency
and simplicity of obtain a solution through a sparse linear system.
Going beyond the present work, our general framework can be adapted
to all previous application areas of the eikonal equation, and, as
alluded to earlier, the variational objective can be readily modified
for to incorporate other constraints that may lead to better reconstructions.

\section{Conclusion}
\label{Section:Discussion}
The Hamilton-Jacobi equation, particularly its specialized form as
the eikonal equation, is at the heart of numerous applications in
vision (shape-from-shading, path planning, medial axis, etc.), and
spurred the rapid development of several innovative computational
techniques to directly solve this \emph{nonlinear} PDE, including
fast marching, and fast sweeping. However, lost in this flurry of
advancing nonlinear solvers was a completely alternative approach,
one which allows you to rigorously approximate solutions to the nonlinear
eikonal as a limiting case of the solution to a corresponding
linear Schr\"odinger equation.  Instead of directly solving the
eikonal equation, the Schr\"odinger formalism results in a generalized,
screened Poisson equation which is solved at very small values of
$\hbar$. In addition, a direct consequence
of our mathematical formulation is that viscosity solutions are naturally
incorporated and obtained when solving the linear differential equation---allowing
one to circumvent explicit viscosity constructions required for any
method that tries to directly solve the nonlinear eikonal.

We initially developed a fast and efficient perturbation series
method for solving the generalized, screened Poisson equation which is guaranteed to converge provided the forcing function
$f$ is positive and bounded. Using the perturbation method 
and the relation (\ref{SfromPhi}), we obtained the solution for the eikonal equation
without spatially discretizing the operators. We later saw that the spatial discretization of the Laplacian operator resulted in a sparse, linear system using which we developed novel solutions to the classical all-pairs, shortest path
problem (a.k.a. path planning). We also illustrated results on shape-from-shading and vessel
centerline extraction. Our approach is straightforward to implement
(by deploying any sparse linear solver) and holds its own against
contemporary fast marching and fast sweeping methods while possessing
the considerable advantage of linearity.  

Our Schr\"odinger-based approach follows the pioneering Hamilton-Jacobi
solvers such as the fast sweeping \cite{Zhao05} and fast marching
\cite{Osher88} methods with the crucial difference being
its linearity.  In future work, we plan to revisit past uses of the eikonal equation and examine
improvements gained by the adoption of our linear framework. We are also
investigating extensions of this approach to other areas such as control
theory. 

\bibliographystyle{amsplain}
\bibliography{Eikonal_JMIVSubmission}

\providecommand{\bysame}{\leavevmode\hbox to3em{\hrulefill}\thinspace}
\providecommand{\MR}{\relax\ifhmode\unskip\space\fi MR }
% \MRhref is called by the amsart/book/proc definition of \MR.
\providecommand{\MRhref}[2]{%
  \href{http://www.ams.org/mathscinet-getitem?mr=#1}{#2}
}
\providecommand{\href}[2]{#2}
\begin{thebibliography}{10}

\bibitem{Abramowitz64}
M.~Abramowitz and I.~A. Stegun, \emph{Handbook of mathematical functions with
  formulas, graphs and mathematical tables}, Dover, New York, NY, 1964.

\bibitem{Arnold89}
V.~I. Arnold, \emph{Mathematical methods of classical mechanics}, Springer, New
  York, NY, 1989.

\bibitem{Basdevant07}
J.~L. Basdevant, \emph{Variational principles in physics}, Springer, New York,
  NY, 2007.

\bibitem{Bracewell99}
R.~N. Bracewell, \emph{The {F}ourier transform and its applications}, 3rd ed.,
  McGraw-Hill, New York, NY, 1999.

\bibitem{Bruss82}
A.~R. Bruss, \emph{The eikonal equation: some results applicable to computer
  vision}, J. Math. Phys. \textbf{23} (1982), no.~5, 890--896.

\bibitem{Butterfield05}
J.~Butterfield, \emph{On {H}amilton-{J}acobi theory as a classical root of
  quantum theory}, Quo-{V}adis {Q}uantum {M}echanics, Springer, New York, NY,
  2005, pp.~239--274.

\bibitem{Canny99}
J.~F. Canny, \emph{Complexity of robot motion planning}, The MIT Press,
  Cambridge, MA, 1988.

\bibitem{Chartier05}
G.~Chartier, \emph{Introduction to optics}, Springer, New York, NY, 2005.

\bibitem{Cooley65}
J.~W. Cooley and J.~W. Tukey, \emph{An algorithm for the machine calculation of
  complex {F}ourier series}, Math. Comp. \textbf{19} (1965), no.~90, 297--301.

\bibitem{Cormen01}
T.~H. Cormen, C.~E. Leiserson, R.~L. Rivest, and C.~Stein, \emph{Introduction
  to algorithms}, 2nd ed., The MIT Press, Cambridge, MA, September 2001.

\bibitem{Crandall92}
M.~G. Crandall, H.~Ishii, and P.L. Lions, \emph{User's guide to viscosity
  solutions of second order partial differential equations}, Bulletin of the
  American Mathematical Society \textbf{27} (1992), no.~1, 1--67.

\bibitem{Demmel97}
J.~W. Demmel, \emph{Applied numerical linear algebra}, SIAM, 1997.

\bibitem{Cohen2001}
T.~Deschamps and L.~D. Cohen, \emph{Fast extraction of minimal paths in {3D}
  images and applications to virtual endoscopy}, Med. Image Anal. \textbf{5}
  (2001), no.~4, 281--299.

\bibitem{Fernandez00}
F.~M. Fernandez, \emph{Introduction to perturbation theory in quantum
  mechanics}, CRC press, 2000.

\bibitem{Fetter03}
A.~L. Fetter and J.~D. Walecka, \emph{Theoretical mechanics of particles and
  continua}, Dover, New York, NY, 2003.

\bibitem{Fousse07}
L.~Fousse, G.~Hanrot, V.~Lef\'evre, P.~P\'elissier, and P.~Zimmermann,
  \emph{{MPFR}: A multiple-precision binary floating-point library with correct
  rounding}, ACM Trans. Math. Softw. \textbf{33} (2007), 1--15.

\bibitem{Goldstein02}
H.~Goldstein, C.P. Poole, and J.~L. Safko, \emph{Classical mechanics}, 3rd ed.,
  Addison Wesley, Boston, MA, 2001.

\bibitem{GMP}
T.~Granlund and \emph{et al.}, \emph{{GNU MP}: {T}he {GNU} {M}ultiple
  {P}recision {A}rithmetic {L}ibrary}, 2012, {url:http://gmplib.org/}.

\bibitem{Griffiths04}
D.~J. Griffiths, \emph{Introduction to quantum mechanics}, 2nd ed., Prentice
  Hall, Upper Saddle River, NJ, 2005.

\bibitem{Gurumoorthy09}
K.~S. Gurumoorthy and A.~Rangarajan, \emph{A {S}chr\"odinger equation for the
  fast computation of approximate {E}uclidean distance functions}, SSVM, LNCS,
  vol. 5567, Springer, 2009, pp.~100--111.

\bibitem{Khatib86}
O.~Khatib, \emph{Real-time obstacle avoidance for manipulators and mobile
  robots}, Int. J. Robot. Res. \textbf{5} (1986), no.~1, 90--98.

\bibitem{Kimmel01}
R.~Kimmel and J.~A. Sethian, \emph{Optimal algorithm for shape from shading and
  path planning}, J. Math. Imaging Vision \textbf{14} (2001), 237--244.

\bibitem{Newton82}
R.G. Newton, \emph{Scattering {T}heory of {W}aves and {P}articles}, 2nd ed.,
  Springer-Verlag, New York, 1982.

\bibitem{Osher02}
S.~J. Osher and R.~P. Fedkiw, \emph{Level set methods and dynamic implicit
  surfaces}, Springer-Verlag, New York, NY, October 2003.

\bibitem{Osher88}
S.~J. Osher and J.~A. Sethian, \emph{Fronts propagating with curvature
  dependent speed: Algorithms based on {H}amilton-{J}acobi formulations}, J.
  Comp. Phys. \textbf{79} (1988), no.~1, 12--49.

\bibitem{Pardos06}
E.~Pardos and O.~Faugeras, \emph{Handbook of mathematical models in computer
  vision}, ch.~Shape from Shading, pp.~275--388, Springer, 2006.

\bibitem{Paris69}
D.T. Paris and F.K. Hurd, \emph{Basic {E}lectromagnetic {T}heory}, McGraw-Hill
  Education, 1969.

\bibitem{Saad03}
Y.~Saad, \emph{Iterative methods for sparse linear systems}, 2nd ed., SIAM,
  2003.

\bibitem{Sethi12}
M.~Sethi, A.~Rangarajan, and K.~S. Gurumoorthy, \emph{The {S}chr\"odinger
  {D}istance {T}ransform {(SDT)} for point-sets and curves}, CVPR, IEEE
  Computer Society, 2012, pp.~198--205.

\bibitem{Sethian96}
J.~A. Sethian, \emph{A fast marching level set method for monotonically
  advancing fronts}, Proc. Nat. Acad. Sci. \textbf{93} (1996), no.~4,
  1591--1595.

\bibitem{Whitham99}
G.~B. Whitham, \emph{Linear and nonlinear waves}, Pure and Applied Mathematics,
  Wiley-Interscience, 1999.

\bibitem{Yatziv06}
L.~Yatziv, A.~Bartesaghi, and G.~Sapiro, \emph{O({N}) implementation of the
  fast marching algorithm}, J. Comp. Phys. \textbf{212} (2006), no.~2,
  393--399.

\bibitem{Zhao05}
H.~K. Zhao, \emph{A fast sweeping method for eikonal equations}, Math. Comp.
  \textbf{74} (2005), no.~250, 603--627.

\end{thebibliography}
\end{document}